%% file: Hitzer_MVDC.tex
\documentclass[twoside]{reporteq}
\usepackage{svcon2e}
\usepackage{makeidx}

\usepackage{latexsym}
\usepackage{amsmath}   
\usepackage{amsfonts}
\usepackage{amssymb}
\usepackage{graphicx}
\usepackage{amsthm}
\usepackage{color}
\input{Makro98}
\input{psfig}
\newcommand{\ed}{\end{document}}
\makeindex
\renewcommand{\theequation}{\arabic{section}.\arabic{equation}}

\newtheorem{prop}{\color{black} Proposition}
\newtheorem{defin}[prop]{\color{black} Definition}
\newtheorem{rem}[prop]{\color{black} Remark}

\newcommand{\bdf}{\begin{defin}}
\newcommand{\edf}{\end{defin}}
\newcommand{\bpro}{\begin{prop}}
\newcommand{\epro}{\end{prop}}
\newcommand{\brem}{\begin{rem}}
\newcommand{\erem}{\end{rem}}

\newcommand{\beqn}{\begin{equation}}
\newcommand{\eeqn}{\end{equation}}
\newcommand{\bea}{\begin{eqnarray}}
\newcommand{\eea}{\end{eqnarray}}
\newcommand{\prf}[1]{\noindent \textbf{\textcolor{black}{Proof} }\textbf{\textcolor{black}{#1}}}

\newcommand{\Slash}[1]{\ensuremath{/\!\!\!\! #1}}

\newcommand{\more}{}
\newcommand{\alles}{}

\newcommand{\cred}{\color{black}}
\definecolor{darkmagenta}{rgb}{0,0,0}
\definecolor{darkblue}{rgb}{0,0,0}
\definecolor{darkgreen}{rgb}{0,0,0}

\begin{document}
\title{Multivector Differential \\ Calculus }
\author{Eckhard M.S. Hitzer\footnote{\color{black}Department of Mechanical Engineering, Fukui University, 
3-9-1 Bunkyo, 910-0024 Fukui, Japan, e-mail: hitzer@mech.fukui-u.ac.jp, homepage: http://sinai.mech.fukui-u.ac.jp/}}
\date{
26 September 2002}
\maketitle
%
%
%







\begin{abstract}
Universal \textcolor{black}{geometric calculus} simplifies and unifies the structure and notation 
of mathematics for all of science and 
engineering, and for technological applications. This paper treats the fundamentals 
of the \textcolor{black}{multivector differential calculus} 
part of geometric calculus.  
The 
\textcolor{black}{multivector differential} is introduced, 
followed by the \textcolor{black}{multivector derivative and the adjoint} of multivector functions. 
The \textcolor{black}{basic rules of multivector differentiation} 
are derived explicitly, as well as a variety of \textcolor{black}{basic multivector derivatives}. 
Finally \textcolor{black}{factorization}, which relates functions of 
vector variables and multivector variables is discussed, and the concepts of both 
\textcolor{black}{simplicial variables and derivatives} are explained.
Everything is proven explicitly in a very elementary level step by step approach. 
The paper is thus intended to serve as 
\textcolor{black}{reference material}, providing a number of details, which are usually skipped in more 
advanced discussions of the subject 
matter. The arrangement of the material closely follows \textcolor{black}{chapter 2 of~\cite{DH:CtoG}}.
\end{abstract}


\begin{quotation}
 But the gate to life is narrow and the way that leads to it is hard,
and there are few people who find it. \ldots I assure you that unless you change
and become like children, you will never enter the Kingdom of heaven.

\emph{\cred Jesus Christ}~\cite{JC:bible}
\end{quotation}

\begin{quotation}
\ldots for geometry, you know, is the gate of science, and the gate
is so low and small that one can only enter it as a little child.

\emph{\cred William K. Clifford}~\cite{Gu:inotr}
\end{quotation}

\section{\cred Introduction}

The work\footnote{First presented in the Geometry session of the 6th Int. Clifford Conference, 
Cookeville, 23 May 2002.}
 is presented as a series of 
{\color{black}definitions (D. or Def.), propositions (P.) and Remarks}. 
\textit{All propositions} (except P. \ref{P67})
\textit{are explicitly proved} in a sometimes
elaborate step by step procedure. 
Such a degree of explicitness seems to be lacking in the literature on the subject\footnote{
Compare e.g. chapter 2 of \cite{DH:CtoG}, the section on multivector calculus of \cite{IB:MVMet} 
or~\cite{DH:MVCalc}.}
 which 
I have come across so far. This motivated me to undertake this study.
It is only assumed that the readers
are familiar with 
{\color{black}geometric (multivector) algebra}, as presented in chapter 1 of  \cite{DH:CtoG}, 
the introductory sections of \cite{EH:VDcalc}, \cite{EH:AMRB} or in numerous other publications on geoemtric 
(Clifford) algebra.
\more

The presentation closely follows the
arrangement in chapter 2 of \cite{DH:CtoG}. It is infact the major purpose of this work to help non-experts
work thoroughly through such a text. There is therefore no great claim to originality, apart from the hope
that it will assist and motivate non-experts and new-comers to delve into the material and become confident
of \textit{comprehensively} understanding and mastering it. The multivector
derivative is to some extent just a generalization of the vector derivative.
I encourage non-experts and new-comers therefore to {\color{black}study \cite{EH:VDcalc}}, 
because it contains a similar 
study {\color{black}restricted to the vector derivative of multivector functions\footnote{
In this context compare also section 2-8 of~\cite{DH:NF1} on directional derivatives,
 section II-2 of~\cite{DH:STC} on vector derivatives and differentials,
 section V-20 of~\cite{DH:STA} on differentiation with its applications throughout~\cite{DH:STA},
 chapter 2 of~\cite{DH:NFII} on geometric calculus, and section 5 of~\cite{GS:SiCalc} on the vector derivative.}}. 
 
In the next section we will start with introducing the {\color{black}multivector differential of multivector functions}, 
i.e. functions with arguments and ranges in the universal geometric algebra ${\cal G}(I)$ with
unit pseudoscalar $I$. 
This will be followed by sections defining the {\color{black}multivector derivative} and the 
{\color{black}adjoint
of multivector functions}. Both the multivector differential and the adjoint can be understood
as two linear approximations associated pointwise to each continuously differential multivector
function.

Conventional scalar differential calculus does not distinguish the three concepts of multivector derivative,
differential and adjoint, because in the scalar case this distinction becomes trivial.

{\color{black}Standard definitions of continuity and scalar differentiability} apply to multivector-valued functions,
because the scalar product determines a \linebreak unique "distance" $|A-B|$ between two elements $A,B\in {\cal G}(I)$. 
\alles

\section{\cred Multivector Differential}

\begin{defin}[\normalfont\textit{differential, $A$-derivative}]
\normalfont
\textup{}For a multivector function $F=F(X)$ defined on the geometric algebra ${\cal G}(I)$:
\begin{equation}
F: X\in {\cal G}(I) \rightarrow F(X) \in {\cal G}(I), I=\langle I \rangle_n,
\end{equation}
\textup{}and for a multivector $A$ and the projection $P(A)=(A\cdot I)\cdot\tilde{I} $ the $A$\textit{-derivative 
(or differential)} is defined by:
\begin{equation}
A\ast\partial_X F(X)\equiv \partial_{\tau} \left. F(X+\tau P(A))\right|_{\tau=0}=
\lim_{\tau \rightarrow 0}\frac{F(X+\tau P(A))-F(X)}{\tau}. 
\end{equation}
\textup{$\ast$ signifies the scalar product defined in \cite{DH:CtoG}, chap. 1, p. 13, (1.44).
By its definition the $A$-derivative is a linear function of $A$ denoted in various ways as:}
\begin{equation}
\underline{F}=\underline{F}(A)=\underline{F}(X,A)=F_A(X)=A\ast\partial_X F(X)=A\ast\partial F.
\end{equation}
\label{D1}
\alles
\end{defin}

\brem
\normalfont
Note the important convention, that inner ($\cdot$), outer ($\wedge$) and scalar ($\ast$) products
have priority over geometric products. 
If $X$ is restricted to a certain subspace of ${\cal G}$, then the projection $P$ in definition \ref{D1}
becomes the projection into that subspace. This is e.g. the case at the beginning of section \ref{sc:BMD}.
\label{R2}
\erem

\begin{prop}
\label{P2}
\begin{equation}
  \underline{F}(A)=\underline{F}(P(A)).
\end{equation}
\end{prop}
\prf{\ref{P2}} 
\begin{align}
  \underline{F}(A)
&\stackrel{\mathrm{Def. \ref{D1}}}{=}\partial_{\tau}\left. F(X+\tau P(A))\right|_{\tau=0}
\nonumber \\
&=
 \partial_{\tau}\left. F(X+\tau P(P(A)))\right|_{\tau=0}\stackrel{\mathrm{Def. \ref{D1}}}{=}\underline{F}(P(A)).
\end{align}

\bpro[{\normalfont {\itshape grade invariance}}]
\label{P3}
\begin{equation}
A\ast\partial\langle F\rangle_r =\langle A\ast\partial F\rangle_r =\langle \underline{F}(A)\rangle_r. 
\end{equation}
\epro
\prf{\ref{P3}} 
\begin{align}
A\ast\partial\langle F\rangle_r  
&\stackrel{\mathrm{Def. \ref{D1}}}{=} 
\partial_{\tau}\langle\left. F(X+\tau P(A))\right\rangle_r|_{\tau=0}
\nonumber \\
&=
\lim_{\tau \rightarrow 0}\frac{\langle F(X+\tau P(A))\rangle_r-\langle F(X)\rangle_r}{\tau} 
\nonumber \\
&\stackrel{\tau\,\,\, \mathrm{ scalar}}{=}
\left\langle \lim_{\tau \rightarrow 0}\frac{F(X+\tau P(A))-F(X)}{\tau} \right\rangle_r
\nonumber \\
&\stackrel{\mathrm{Def. \ref{D1}}}{=}
\langle A\ast\partial F\rangle_r 
\stackrel{\mathrm{Def. \ref{D1}}}{=}
\langle \underline{F}(A)\rangle_r. 
\end{align}

\bpro
\begin{equation}
\underline{F}(A+B)=\underline{F}(A)+\underline{F}(B).
\end{equation}
\label{P4}
\epro
\prf{\ref{P4}}
\begin{align}
  \underline{F}(A+B)
&\stackrel{\mathrm{Def. \ref{D1} }}{=} 
\lim_{\tau \rightarrow 0}\frac{F(X+\tau P(A+B))-F(X)}{\tau}
\nonumber \\
&\hspace{-0.5cm}\stackrel{\mbox{\scriptsize linearity of $P$}}{=}
\lim_{\tau \rightarrow 0}\frac{F(X+\tau P(A)+\tau P(B))-F(X)}{\tau} 
=\lim_{\tau \rightarrow 0}
\nonumber \\
&\hspace{-1.7cm}\left\{ \frac{F(X+\tau P(A)+\tau P(B))-F(X+\tau P(B))}{\tau} 
    +\frac{F(X+\tau P(B))-F(X)}{\tau} \right\} 
\nonumber \\
&\hspace{-1cm}\stackrel{\mbox{\scriptsize $X_{\tau}\equiv X+\tau P(B),$ Def. \ref{D1}}}{=}
\lim_{\tau \rightarrow 0}\frac{F(X_{\tau}+\tau P(A))-F(X_{\tau})}{\tau} +\underline{F}(B)
\nonumber \\
&\hspace{-1cm}\stackrel{\mbox{\scriptsize 
                                $\lim_{\tau \rightarrow 0}X_{\tau}=X,$ Def. \ref{D1} }}{=}
A\ast\partial_X F(X)+\underline{F}(B) 
\stackrel{\mathrm{Def. \ref{D1}}}{=}\underline{F}(A)+\underline{F}(B).
\end{align}
An alternative rigorous proof of P. \ref{P4} is the following: I first define a function $\varepsilon_F(X,A,\tau)$,
which is according to Def. \ref{D1} continuous at $\tau=0$
\begin{equation}
\varepsilon_F(X,A,\tau) \equiv 
\left\{
\begin{array}{l}
\frac{F(X+\tau P(A))-F(X)}{\tau}-\underline{F}(X,A)\,\,\,\,\, \tau\neq 0, \\
0 \,\,\,\,\, \tau=0.
\end{array}
\right.
\label{eq:P4.A1}
\end{equation}
We hence have
\begin{align}
F(X+\tau P(A)) 
\stackrel{\mbox{\scriptsize (\ref{eq:P4.A1})  }}{=}
F(X)+\tau \underline{F}(X,A) + \tau \varepsilon_F(X,A,\tau) 
\label{eq:P4.A2}
\end{align}
and
\begin{align}
&F(X+\tau P(A+B))
\stackrel{\mbox{\scriptsize linearity of $P$  }}{=}
F(\underbrace{X+\tau P(A)} + \tau P(B)) 
\nonumber \\
&\stackrel{\mbox{\scriptsize (\ref{eq:P4.A2})  }}{=}
F(X+\tau P(A)) + \tau \underline{F}(X+\tau P(A),B) +\tau \varepsilon_F(X+\tau P(A),B,\tau)
\nonumber \\
&\stackrel{\mbox{\scriptsize (\ref{eq:P4.A2})  }}{=}
  F(X)+\tau \underline{F}(X,A) + \tau \varepsilon_F(X,A,\tau)
  \nonumber \\
  &\hspace{1cm}+\tau \underline{F}(X+\tau P(A),B) +\tau \varepsilon_F(X+\tau P(A),B,\tau).
\label{eq:P4.A3}
\end{align}
We can now calculate $\underline{F}(X,A+B)$ according to Def. \ref{D1}
\begin{align}
&\underline{F}(X,A+B)
\stackrel{\mbox{\scriptsize Def. \ref{D1}  }}{=}
\lim_{\tau \rightarrow 0}
  \frac{F(X +\tau P(A+B))-F(X)}{\tau}
\stackrel{\mbox{\scriptsize (\ref{eq:P4.A3})  }}{=}
\lim_{\tau \rightarrow 0}
\nonumber \\
&
 \{\underline{F}(X,A) + \varepsilon_F(X,A,\tau)
+\underline{F}(X+\tau P(A),B) +\varepsilon_F(X+\tau P(A),B,\tau)\}
\nonumber \\
&=\lim_{\tau \rightarrow 0}
\{\underline{F}(X,A) + \varepsilon_F(X,A,\tau)
+\underline{F}(X,B) +\varepsilon_F(X,B,\tau)\}
\nonumber \\
&\stackrel{\mbox{\scriptsize (\ref{eq:P4.A1})  }}{=}
\underline{F}(X,A)+\underline{F}(X,B).
\end{align}

\bpro 
\normalfont
For constant scalar $\lambda$
\begin{equation}
\underline{F}(\lambda A)=\lambda \underline{F}(A).
\end{equation}
\label{P5}
\epro
\prf{\ref{P5}}
\begin{align}
\hspace{-2.7cm}\underline{F}(\lambda A)
&\stackrel{\mathrm{Def. \ref{D1}}}{=}
\lim_{\tau \rightarrow 0}\frac{F(X+ \tau \lambda P(A))-F(X)}{\tau}
\nonumber \\
&\stackrel{\lambda \neq 0}{=}
\lim_{\tau \rightarrow 0}\lambda \frac{F(X+ \tau \lambda P(A))-F(X)}{\tau\lambda} 
\nonumber \\
&\hspace{-0.3cm}\stackrel{\tau \rightarrow {\tau\prime}= \tau\lambda}{=}
\lambda \lim_{{\tau\prime} \rightarrow 0}\frac{F(X+ {\tau\prime} P(A))-F(X)}{{\tau\prime}}
\stackrel{\mathrm{Def. \ref{D1}}}{=}
\lambda \underline{F}(A),
\end{align}
in case of $\lambda =0$ we simply have
\begin{equation}
\underline{F}(0 A)=\underline{F}(0)\stackrel{\mathrm{Def. \ref{D1}}}{=}
\lim_{\tau \rightarrow 0}\frac{F(X)-F(X)}{\tau}=0=0\underline{F}(A).
\end{equation}

\bpro[\normalfont \textit{sum rule}]
\normalfont
For two multivector functions $F,G$ on ${\cal G}(I)$:
\begin{equation}
  A\ast\partial_X(F+G)= A\ast\partial_XF+ A\ast\partial_XG.
\end{equation}
\label{P6}
\epro
\prf{\ref{P6}}
\bea
  \lefteqn{A\ast\partial_X(F+G) \stackrel{\mathrm{Def. \ref{D1}}}{=}
\partial_{\tau}[F(X+ \tau P(A))+G(X+ \tau P(A))]=} \\
&&=\partial_{\tau}F(X+ \tau P(A))+\partial_{\tau}G(X+ \tau P(A)) \stackrel{\mathrm{Def. \ref{D1}}}{=}
 A \ast \partial_X F+ A\ast\partial_X G. \nonumber
\eea

\bpro[{\normalfont {\itshape product rule}}]
\normalfont
For two multivector functions $F,G$ on ${\cal G}(I)$:
\begin{equation}
A\ast\partial_X(FG)=(A\ast\partial_X F)G+F(A\ast\partial_X G).
\end{equation}
\label{P7}
\epro
\prf{\ref{P7}}
\bea
\lefteqn{ A\ast\partial_X(FG)  \stackrel{\mathrm{Def. \ref{D1}}}{=} 
\lim_{\tau \rightarrow 0}
\frac{F(X+\tau P(A))G(X+\tau P(A))-F(X)G(X)}{\tau}}
\nonumber \\
&&\hspace{7mm}= \lim_{\tau \rightarrow 0}  \frac{F(X+\tau P(A))G(X+\tau P(A))-F(X)G(X+\tau P(A))}{}\cdots
\nonumber \\ &&\hspace{17mm}
 \frac{+F(X)G(X+\tau P(A))-F(X)G(X)}{\tau}  \nonumber \\
&&\hspace{7mm}=\lim_{\tau \rightarrow 0} \left\{ \frac{F(X+\tau P(A))-F(X)}{\tau}\,\,G(X+\tau P(A))+
\right.
\nonumber \\
&&\hspace{53mm}\left.+F(X) \,\frac{G(X+\tau P(A))-G(X)}{\tau} \right\}  \nonumber \\
&&\hspace{7mm}\stackrel{\mathrm{Def. \ref{D1}}}{=}
(A\ast\partial_X F)\lim_{\tau \rightarrow 0}G(X+\tau P(A))+F\,(A\ast\partial_X G)  \nonumber \\
&&\hspace{7mm}=(A\ast\partial_X F)G+F\,(A\ast\partial_X G).
\eea
\alles

\bpro[{\normalfont {\itshape constant function}}]
\normalfont
For $B$ independent of $X$:
\begin{equation}
A \ast \partial_X B=0.
\end{equation}
\label{P17}
\epro
\prf{\ref{P17}}
Def. \ref{D1} gives for $F\equiv B$:
\begin{equation}
 A\ast \partial_X B=\lim_{\tau \rightarrow 0} \frac{B-B}{\tau}=0.
\end{equation}

\bpro[{\normalfont {\itshape identity}}]
\begin{equation}
A\ast \partial_XX = P(A).
\end{equation}
\label{P18}
\epro
\prf{\ref{P18}}
For $F(X)\equiv X$:
\begin{equation}
A\ast\partial_XX\stackrel{\mbox{\scriptsize Def. \ref{D1}}}{=}
\lim_{\tau \rightarrow 0} \frac{X+\tau P(A)-X}{\tau}=P(A).
\end{equation}

\bpro
\normalfont
For $B$ independent of $X$:
\begin{equation}
A\ast\partial_X(X\ast B)=P(A)\ast B.
\end{equation}
\label{P19}
\epro
\prf{\ref{P19}}
\begin{align}
A \ast \partial_X(X\ast B) \stackrel{\mbox{\scriptsize Def. of scalar prod.}}{=}
&A \ast \partial_X \langle XB \rangle_0 \stackrel{\mbox{\scriptsize P. \ref{P3}}}{=} 
\langle A \ast \partial_X XB \rangle_0 \stackrel{\mbox{\scriptsize P. \ref{P7}}}{=} \\
\langle (A \ast \partial_X X)B +XA \ast \partial_XB \rangle_0
&\stackrel{\mbox{\scriptsize P. \ref{P17}, P. \ref{P18}}}{=} 
\langle (P(A)B +0 \rangle_0= P(A)\ast B.\nonumber
\end{align}

\bpro[{\normalfont {\itshape Taylor expansion}}]
\begin{equation}
F(X+P(A))=\sum^{\infty}_{k=0}\frac{1}{k!}(A\ast\partial_X)^k F=\exp(A\ast\partial_X)F.
\end{equation}
\label{P8}
\epro
\prf{\ref{P8}}
\begin{equation}
G(\tau)\equiv F(X+\tau P(A)) \Rightarrow \frac{dG(0)}{d\tau}=\partial_{\tau}
\left.F(X+\tau P(A))\right|_{\tau=0}\stackrel{\mathrm{Def. \ref{D1}}}{=} A*\partial_XF, 
\end{equation}
\bea
\lefteqn{\hspace{-1.5mm}\mbox{Def. of }G \Rightarrow \frac{d^2G(0)}{d\tau^2}=\left.\partial_{\tau}\left(\partial_{\tau}
F(X+\tau P(A)) \right)\right|_{\tau=0}} \\
&&\hspace{21mm}\stackrel{\mathrm{Def. \ref{D1}}}{=}
 A*\partial_X \left.\left(\partial_{\tau}F(X+\tau P(A)) \right)\right|_{\tau=0}
\stackrel{\mathrm{Def. \ref{D1}}}{=}
(A*\partial_X)^2F.
\nonumber
\eea
\begin{equation} 
\mbox{In general: }  \frac{d^kG(0)}{d\tau^k}=(A*\partial_X)^kF.
\end{equation}
\more

Taylor expansion of $G$:
\begin{equation}
G(1)=G(0+1)=G(0)+\frac{dG(0)}{d\tau}+\frac{1}{2}\frac{d^2G(0)}{d\tau^2}+\ldots=\sum^{\infty}_{k=0}\frac{1}{k!}
\frac{d^kG(0)}{d\tau^k}.
\end{equation}
\begin{equation}
\Rightarrow G(1)=F(X+P(A))=\sum^{\infty}_{k=0}\frac{1}{k!}(A\ast\partial_X)^k F=\exp(A\ast\partial_X)F.
\end{equation}
\alles

\bpro[{\normalfont {\itshape chain rule}}]
\normalfont
\begin{equation}
 f: X\in {\cal G}(I)\rightarrow X^{\prime} \in {\cal G}(I^{\prime}) \mbox{ with } I^{\prime} 
=\langle I^{\prime} \rangle_n
\end{equation}
The composite function $F(X)=G(f(X))$ has the differential
\begin{equation}
A\ast\partial F=A\ast\partial G(f(X)) =\underline{f}(A)\ast\partial G
\end{equation}
\label{P9}
\epro
\prf{\ref{P9}}
The Taylor expansion (P.\ref{P8}) of $f$ gives:

\begin{align}
f(X+\tau P(A))&=f(X)+(\tau A\ast \partial_X )f+\frac{1}{2}(\tau A\ast \partial_X )^2f+\ldots 
\nonumber   \\
&\hspace{-0.1cm}
\stackrel{\mathrm{Def. \ref{D1}}}{=}
f(X)+\underline{f}(\tau A)+\frac{1}{2}(\tau A\ast \partial_X )^2f+\ldots 
\nonumber \\
&\stackrel{\mathrm{P. \ref{P5}}}{=} f(X)+\tau\underline{f}( A)+\frac{\tau^2}{2}( A\ast \partial_X )^2f+\ldots
\end{align}
Therefore
\begin{align}
A\ast \partial_X F
&=A\ast \partial_X G(f(X))\stackrel{\mathrm{Def. \ref{D1}}}{=}
\partial_{\tau}\left.G(f(X+\tau P(A)))\right|_{\tau=0}
\nonumber \\
&\hspace{-0.2cm}\stackrel{\mathrm{Taylor}}{=}\partial_{\tau}\left.G(f(X)+\tau\underline{f}( A))\right|_{\tau=0}
\nonumber \\
&\hspace{-1.1cm}\stackrel{ \mbox{\scriptsize range of $f,$ $\underline{f}$ in ${\cal G}(I^{\prime})$}}{=}
\partial_{\tau}\left.G(f(X)+\tau P(\underline{f}( A)))\right|_{\tau=0} \nonumber \\
&\hspace{-0.1cm}\stackrel{\mathrm{Def. \ref{D1}}}{=}
\underline{f}(A)\ast \partial_{X^{\prime}} \left. G(X^{\prime}) \right|_{X^{\prime}=f(X)}.
\end{align}
\alles

\brem
\label{R10}
\normalfont
Proof {\ref{P9}} employs the Taylor expansion (P. \ref{P8}). Yet there are continuously differentiable functions
without a Taylor expansion. But the chain rule P. \ref{P9} will still apply as long as a function $f$ has the linear
approximation $f(X+\tau P(A))\approx f(X)+\tau\underline{f}( A)$ for sufficiently small values of $\tau$.\footnote{
I owe this remark to Dr. H. Ishi, of Yokohama City University.} 

I therefore want to give another proof for the 
chain rule based on the linearization of multivector functions with the help of the differential. 
One such linearization is already given in (\ref{eq:P4.A1}) and (\ref{eq:P4.A2}). In the same way I define
\begin{equation}
\varepsilon_f(X,A,\tau) \equiv 
\left\{
\begin{array}{l}
\frac{f(X+\tau P(A))-f(X)}{\tau}-\underline{f}(X,A)\,\,\,\,\, \tau\neq 0, \\
0 \,\,\,\,\, \tau=0,
\end{array}
\right.
\label{eq:R10.1}
\end{equation}
Hence we have the linearization of $f:X\in{\cal G}\rightarrow f(X)\in {\cal G}$
\begin{align}
f(X+\tau P(A)) 
\stackrel{\mbox{\scriptsize (\ref{eq:R10.1})  }}{=}
f(X)+\tau \underline{f}(X,A) + \tau \varepsilon_f(X,A,\tau). 
\label{eq:R10.2}
\end{align}
I further define
\begin{equation}
\varepsilon_G(Y,B,\tau) \equiv 
\left\{
\begin{array}{l}
\frac{G(Y+\tau P(B))-G(Y)}{\tau}-\underline{G}(Y,B)\,\,\,\,\, \tau\neq 0, \\
0 \,\,\,\,\, \tau=0,
\end{array}
\right.
\label{eq:R10.3}
\end{equation}
and hence
\begin{align}
G(Y+\tau P(B)) 
\stackrel{\mbox{\scriptsize (\ref{eq:R10.3})  }}{=}
G(Y)+\tau \underline{G}(Y,B) + \tau \varepsilon_G(Y,B,\tau). 
\label{eq:R10.4}
\end{align}
According to Def. \ref{D1} both functions $\varepsilon_f$ and $\varepsilon_G$ are continuous at 
$\tau=0.$ Considering that 
$f(X)\in {\cal G}\stackrel{\mbox{\scriptsize Def. \ref{D1}}}{\Rightarrow}\underline{f}\in {\cal G}
\stackrel{\mbox{\scriptsize (\ref{eq:R10.1})}}{\Rightarrow}\varepsilon_f\in {\cal G}$ (the same
is valid for a function $f$ on a linear subspace of ${\cal G}$),
we can now rewrite $F(X)=G(f(X))$ as
\begin{align}
&G(f(X+\tau P(A)))
\stackrel{\mbox{\scriptsize (\ref{eq:R10.2})  }}{=}
G(f(X)+\tau \{\underline{f}(X,A) + \varepsilon_f(X,A,\tau)\})
\nonumber \\
&\stackrel{\mbox{\scriptsize $f,\varepsilon_f\in {\cal G}$  }}{=}
G(f(X)+\tau P(\underline{f}(X,A) + \varepsilon_f(X,A,\tau)) )
\nonumber \\
&\stackrel{\mbox{\scriptsize (\ref{eq:R10.4})  }}{=}
G(Y=f(X))+\tau \underline{G}(Y=f(X),B=\underline{f}(X,A) + \varepsilon_f(X,A,\tau))
\nonumber \\
&\hspace{1cm}+\tau \varepsilon_G(f(x),\underline{f}(X,A) + \varepsilon_f(X,A,\tau),\tau). 
\end{align}
Subtracting $G(f(X))$, deviding by $\tau$ and taking the $\lim_{\tau \rightarrow 0}$ we get 
according to Def. \ref{D1} the
differential of the composite function $F(X)=G(f(X))$
\begin{align}
&\lim_{\tau \rightarrow 0} \frac{G(f(X+\tau P(A)))-G(f(X))}{\tau}
=\lim_{\tau \rightarrow 0} 
\nonumber \\
&\{\underline{G}\left(f(X),\underline{f}(X,A) + \varepsilon_f(X,A,\tau)\right)
+\varepsilon_G(f(x),\underline{f}(X,A) + \varepsilon_f(X,A,\tau),\tau)\}
\nonumber \\
&\stackrel{\mbox{\scriptsize (\ref{eq:R10.1}),(\ref{eq:R10.3})  }}{=}
\underline{G}(f(X),\underline{f}(X,A)).
\end{align}
This concludes the alternative proof for the differential of composite functions (chain rule). 
\erem

\bdf[{\normalfont {\itshape second differential}}]
\normalfont
\begin{equation}
F_{AB}=F_{AB}(X)\equiv B\ast \dot{\partial} A\ast \partial \dot{F}.
\end{equation}
The over-dot notation indicates, that $\dot{\partial}$ acts only on $F$.
\label{D10}
\edf

\bpro[{\normalfont {\itshape integrability}}]
\normalfont
Let $F(X+\tau P(A)+\sigma P(B))$ be a twice continously differentiable function in the vicinity of $X,$ i.e.
for all values of $0\leq \tau \leq \tau_1$ and $0 \leq \sigma \leq \sigma_1.$
Then
\begin{equation}
 F_{AB}(X)=F_{BA}(X).
\end{equation}
\label{P11}
\epro

\brem
\normalfont
To motivate the slightly more elaborate proof of P. \ref{P11} which follows, I present this "handwaving" argument:
\bea
\lefteqn{F_{AB}(X)\stackrel{\mathrm{Def. \ref{D10}}}{=}B\ast \dot{\partial} A\ast \partial \dot{F}
\stackrel{\mathrm{Def. \ref{D1}}}{=}B\ast  \dot{\partial}\left(\partial_{\tau}\left.\dot{F}(X+ \tau P(A))\right|_{\tau=0} \right)}
\nonumber \\
&&\hspace{-0.75cm}\stackrel{\mathrm{Def. \ref{D1}}}{=}
\partial_{\sigma}\partial_{\tau}\left.F(X+ \tau P(A)+\sigma P(B))\right|_{\tau=\sigma=0} 
\label{eq:R11.1}
\\
&&\hspace{-0.6cm}=\lim_{\sigma \rightarrow 0}\lim_{\tau \rightarrow 0}
\frac{\frac{F(X+ \tau P(A)+\sigma P(B))-F(X+ \tau P(B))}{\tau}-
\frac{F(X+ \tau P(A)-F(x))}{\tau}}{\sigma} =\lim_{\sigma \rightarrow 0}
\nonumber \\
&&\hspace{-0.8cm}\lim_{\tau \rightarrow 0}
\frac{F(X+ \tau P(A)+\sigma P(B))-F(X+ \tau P(B))-F(X+ \tau P(A))+F(x)}{\sigma \tau}. \nonumber
\eea
The last expression is symmetric under\footnote{
Remark \ref{R17} is "handwaving", because the implied interchange of the limits in 
$\sigma$ and $\tau$ is mathematically nontrivial. I am very much indebted to Dr. H. Ishi (Yokohama City
University, Japan) for discussing P. \ref{P11} with me and for suggesting the proof of it given below.}
 the exchange $(P(A),\tau)\leftrightarrow (P(B),\sigma)$. Hence
\begin{equation}
F_{AB}(X)=B\ast \dot{\partial} A\ast \partial \dot{F}(X)
=A\ast \dot{\partial} B\ast \partial \dot{F}(X)= F_{BA}(X).
\end{equation}
\label{R17}
\erem

\prf{\ref{P11}}
According to the {\itshape fundamental theorem of calculus} (\cite{EDM2}, 216.C and~\cite{OF:ANAL1}, p. 140, Satz 3) we have for in the intervall
$[0,\tau_1]$ continously differentiable functions $g(\tau)$
\begin{equation}
g(\tau_1)-g(0) = \int_0^{\tau_1}\partial_{\tau}g({\eta}) d\eta.
\label{eq:P11.1}
\end{equation}
We now use the in $\tau$ and $\sigma$ twice continuously differentiable function 
\begin{equation}
f(\tau,\sigma)\equiv F(X+\tau P(A)+\sigma P(B))
\label{eq:P11.2}
\end{equation}
to define
\begin{equation}
M(\tau_1,\sigma_1)\equiv f(\tau_1,\sigma_1)-f(\tau_1,0)-f(0,\sigma_1)+f(0,0).
\end{equation}
Following the notation of (\ref{eq:R11.1}) we can therefore write
\begin{align}
&F_{AB}=\lim_{\sigma_1\rightarrow 0}\left( \lim_{\tau_1 \rightarrow 0} 
  \frac{M(\tau_1,\sigma_1)}{\sigma_1\tau_1}\right)
=\left. \partial_{\sigma}\partial_{\tau}f\right|_{\tau=\sigma=0}=\partial_{\sigma}\partial_{\tau}f(0,0),
\label{eq:P11.3}
\\
&F_{BA}=\lim_{\tau_1\rightarrow 0}\left( \lim_{\sigma_1 \rightarrow 0} 
  \frac{M(\tau_1,\sigma_1)}{\sigma_1\tau_1}\right)
=\left. \partial_{\tau}\partial_{\sigma}f\right|_{\tau=\sigma=0}=\partial_{\tau}\partial_{\sigma}f(0,0).
\label{eq:P11.4}
\end{align}
Using the fundamental theorem of calculus twice the function $M(\tau_1,\sigma_1)$ can be reexpressed 
on the one hand by
\begin{align}
&M(\tau_1,\sigma_1)\stackrel{\mbox{\scriptsize (\ref{eq:P11.1})}}{=}
\int_{0}^{\sigma_1}\partial_{\sigma}f(\tau_1,\xi) d\xi - \int_{0}^{\sigma_1}\partial_{\sigma}f(0,\xi) d\xi
\nonumber
\\
&= \int_{0}^{\sigma_1}\left[ \partial_{\sigma}f(\tau_1,\xi) - \partial_{\sigma}f(0,\xi) \right] d\xi
\stackrel{\mbox{\scriptsize (\ref{eq:P11.1})}}{=}
\int_{0}^{\sigma_1}\int_{0}^{\tau_1} \partial_{\tau}\partial_{\sigma}f(\eta,\xi) d\eta d\xi, 
\label{eq:P11.5}
\end{align}
and on the other hand by
\begin{align}
&M(\tau_1,\sigma_1)
=f(\tau_1,\sigma_1)-f(0,\sigma_1)-f(\tau_1,0)+f(0,0)
\nonumber 
\\
&\stackrel{\mbox{\scriptsize (\ref{eq:P11.1})}}{=}
\int_{0}^{\tau_1}\partial_{\tau}f(\eta,\sigma_1) d\eta - \int_{0}^{\tau_1}\partial_{\tau}f(\eta,0)  d\eta
= \int_{0}^{\tau_1}\left[ \partial_{\tau}f(\eta,\sigma_1) - \partial_{\tau}f(\eta,0) \right] d\eta
\nonumber
\\
&\stackrel{\mbox{\scriptsize (\ref{eq:P11.1})}}{=}
\int_{0}^{\tau_1}\int_{0}^{\sigma_1} \partial_{\sigma}\partial_{\tau}f(\eta,\xi) d\xi d\eta. 
\label{eq:P11.6}
\end{align}
It follows from (\ref{eq:P11.5}) and (\ref{eq:P11.6}) that
\begin{align}
M(\tau_1,\sigma_1)
\stackrel{\mbox{\scriptsize (\ref{eq:P11.5})}}{=}
\int_{0}^{\sigma_1}\int_{0}^{\tau_1} \partial_{\tau}\partial_{\sigma}f(\eta,\xi) d\eta d\xi
\stackrel{\mbox{\scriptsize (\ref{eq:P11.6})}}{=}
\int_{0}^{\tau_1}\int_{0}^{\sigma_1} \partial_{\sigma}\partial_{\tau}f(\eta,\xi) d\xi d\eta.
\label{eq:P11.7}
\end{align}
The assumption of the double continous differentiability of $f$ means that 
\begin{align}
\partial_{\tau}\partial_{\sigma}f(\eta,\xi) &\,\,\,\mbox{and}\,\,\, \partial_{\sigma}\partial_{\tau}f(\eta,\xi)
\,\,\,\mbox{  are both continuous for   }\,\,\,
\nonumber \\ 
&0\leq \eta \leq \tau_1, 0\leq \xi \leq \sigma_1.
\end{align}
This can in turn has the consequence that
\begin{align}
m(\tau_1,\sigma_1)\equiv \,\,\,&\mbox{max}_{\{0\leq \eta \leq \tau_1, 0\leq \xi \leq \sigma_1\}}
\left| \partial_{\sigma}\partial_{\tau}f(\eta,\xi)-\partial_{\sigma}\partial_{\tau}f(0,0)\right|
\rightarrow 0
\nonumber \\
&\mbox{for} \,\,\,(\tau_1,\sigma_1) \rightarrow (0,0),
\label{eq:P11.8}
\end{align}
and
\begin{align}
n(\tau_1,\sigma_1)\equiv \,\,\,&\mbox{max}_{\{0\leq \eta \leq \tau_1, 0\leq \xi \leq \sigma_1\}}
\left| \partial_{\tau}\partial_{\sigma}f(\eta,\xi)-\partial_{\tau}\partial_{\sigma}f(0,0)\right|
\rightarrow 0
\nonumber \\
&\mbox{for} \,\,\,(\tau_1,\sigma_1) \rightarrow (0,0).
\label{eq:P11.9}
\end{align}
We therefore have
\begin{align}
&\left| \frac{M(\tau_1,\sigma_1)}{\sigma_1\tau_1}-\partial_{\tau}\partial_{\sigma}f(0,0) \right|
\nonumber \\
&\stackrel{\mbox{\scriptsize (\ref{eq:P11.5})}}{=}
\left| 
 \frac{1}{\sigma_1\tau_1}\int_{0}^{\sigma_1}\int_{0}^{\tau_1} \partial_{\tau}\partial_{\sigma}f(\eta,\xi) d\eta d\xi
-\frac{1}{\sigma_1\tau_1}\int_{0}^{\sigma_1}\int_{0}^{\tau_1} \partial_{\tau}\partial_{\sigma}f(0,0) d\eta d\xi
\right|
\nonumber \\
&\hspace{2mm}\leq 
\frac{1}{\sigma_1\tau_1}\int_{0}^{\sigma_1}\int_{0}^{\tau_1}\left| 
\partial_{\tau}\partial_{\sigma}f(\eta,\xi)-\partial_{\tau}\partial_{\sigma}f(0,0)
\right| d\eta d\xi
\nonumber \\
&\stackrel{\mbox{\scriptsize (\ref{eq:P11.9})}}{\leq}
\frac{1}{\sigma_1\tau_1}\int_{0}^{\sigma_1}\int_{0}^{\tau_1} n(\tau_1,\sigma_1) d\eta d\xi
\nonumber \\
&\hspace{2mm}=\,\,\,
n(\tau_1,\sigma_1) \frac{\int_{0}^{\sigma_1}\int_{0}^{\tau_1}d\eta d\xi}{\sigma_1\tau_1}
=
n(\tau_1,\sigma_1).
\end{align}
Because of (\ref{eq:P11.9}) we then get
\begin{align}
\lim_{(\tau_1,\sigma_1)\rightarrow (0,0)}
\left| \frac{M(\tau_1,\sigma_1)}{\sigma_1\tau_1}-\partial_{\tau}\partial_{\sigma}f(0,0)\right|
\stackrel{\mbox{\scriptsize (\ref{eq:P11.9})}}{=}0.
\label{eq:P11.10}
\end{align}
Likewise we have
\begin{align}
&\left| \frac{M(\tau_1,\sigma_1)}{\sigma_1\tau_1}-\partial_{\sigma}\partial_{\tau}f(0,0)\right|
\nonumber \\
&\stackrel{\mbox{\scriptsize (\ref{eq:P11.6})}}{=}
\left|
 \frac{1}{\sigma_1\tau_1}\int_{0}^{\tau_1}\int_{0}^{\sigma_1} \partial_{\sigma}\partial_{\tau}f(\eta,\xi) d\xi d\eta
-\frac{1}{\sigma_1\tau_1}\int_{0}^{\tau_1}\int_{0}^{\sigma_1} \partial_{\sigma}\partial_{\tau}f(0,0) d\xi d\eta
\right|
\nonumber \\
&\hspace{2mm}\leq 
\frac{1}{\sigma_1\tau_1}\int_{0}^{\tau_1}\int_{0}^{\sigma_1}\left| 
\partial_{\sigma}\partial_{\tau}f(\eta,\xi)-\partial_{\sigma}\partial_{\tau}f(0,0)
\right| d\xi d\eta
\nonumber \\
&\stackrel{\mbox{\scriptsize (\ref{eq:P11.8})}}{\leq}
\frac{1}{\sigma_1\tau_1}\int_{0}^{\tau_1}\int_{0}^{\sigma_1} m(\tau_1,\sigma_1) d\xi d\eta
\nonumber \\
&\hspace{2mm}=\,\,\,
m(\tau_1,\sigma_1) \frac{\int_{0}^{\tau_1}\int_{0}^{\sigma_1}d\xi d\eta}{\sigma_1\tau_1}
=
m(\tau_1,\sigma_1).
\end{align}
Because of (\ref{eq:P11.8}) we then get
\begin{align}
\lim_{(\tau_1,\sigma_1)\rightarrow (0,0)}
\left| \frac{M(\tau_1,\sigma_1)}{\sigma_1\tau_1}-\partial_{\sigma}\partial_{\tau}f(0,0)\right|
\stackrel{\mbox{\scriptsize (\ref{eq:P11.8})}}{=}0.
\label{eq:P11.11}
\end{align}
Equations (\ref{eq:P11.7}), (\ref{eq:P11.10}) and (\ref{eq:P11.11}) show that
\begin{align}
\partial_{\tau}\partial_{\sigma}f(0,0)
=
\partial_{\sigma}\partial_{\tau}f(0,0).
\label{eq:P11.12}
\end{align}
Finally (\ref{eq:P11.12}) results with (\ref{eq:P11.3}) and (\ref{eq:P11.4}) in the proof of P. \ref{P11}:
\begin{align}
F_{AB}(X) = F_{BA}(X).
\end{align}

\section{\cred Multivector Derivative}

\bdf
\normalfont
The brackets $\langle \ldots \rangle_r$ indicate the selection of the $r$-grade part of the multivector 
expression enclosed by them.
\begin{equation}
A_{\bar{r}}\equiv \langle A \rangle_{r}, \hspace{1cm} \langle A \rangle\equiv \langle A \rangle_0=A_{\bar{0}}.
\end{equation}
\label{D15}
\edf

\bdf[{\normalfont {\itshape derivative}}]
\normalfont
The \textit{derivative}\footnote{
It is important to note that based on directed integration, which uses a "directed Riemann measure", one can define
the multivector derivative by a limit process applied to an integral over the boundary of
 a smooth open $r$-dimensional surface "tangent" to (in this definition $r$-vector) $X$. 
In the limit process the volume of the $r$-dimensional surface is shrunk to zero. 
(Compare~\cite{DH:MVCalc}, section 4.) For an analogous definition of the vector derivative see section 
5 of~\cite{GS:SiCalc}.
} of a multivector function $F(X)$ by its argument $X$
\begin{equation}
\partial_X F(X)=\partial F,
\end{equation}
with derivative operator $\partial_X$, is assumed
\begin{description}
 \item[(i)]  to have the algebraic properties of a multivector in ${\cal G}(I)$, with $I$ the unit pseudoscalar.
 \item[(ii)] that as in definition \ref{D1}: $A\ast \partial_X$ with $A\in {\cal G}(I)$ equals the differential
\begin{equation}
  A\ast \partial_X F = \underline{F}_A=\sum_{r=0}^{n}A_{\bar{r}}*\partial_X F
=\sum_{r=0}^{n}A*\partial_{\bar{r}} F=\sum_{r=0}^{n}A_{\bar{r}}*\partial_{\bar{r}} F.
\end{equation}
\normalfont [Compare \cite{DH:CtoG}, p. 13, (1.46) and (1.47a).]
\end{description}
\label{D12}
\edf

\brem
\normalfont
Property (ii) in Def. \ref{D12} expresses that the derivative operator acts like a map from the space 
(usually ${\cal G}$ or a linear subspace of ${\cal G}$) of multivectors $A$ to the space of 
(differentials of) multivector functions $F$:
\begin{align}
\partial_X F: A\in {\cal G} \rightarrow A\ast\partial_X F\in {\cal G}.
\end{align}
In order to avoid easily occuring confusions, I want to point out that Remark \ref{R2} clearly implies that
\begin{align}
A \ast \partial_X F &= (A\ast \partial_X)F
\nonumber \\
& \neq A \ast (\partial_X F).
\end{align}
\label{R13}
\erem

\bpro[{\normalfont {\itshape algebraic properties of $\partial_X$}}]
\normalfont
\begin{equation}
\partial_X=P(\partial_X)=\sum_J a^J a_J\ast\partial_X
\stackrel{\mbox{\scriptsize for $a^J=const.$}}{=}
\sum_J a_J \ast \partial_X \,\,a^J,
\label{eq:P13.1a}
\end{equation}
with $P$ the multivector projection into ${\cal G}(\vec{a}_1\wedge \vec{a}_2\wedge\ldots \wedge \vec{a}_n)$.
$a_J$ is a simple blade basis\footnote{The vectors $\vec{a}_{j_k}, k=1,\ldots, n$ 
in (\ref{eq:P13.2}) are not necessarily orthonormal, but can be orthogonalized by a procedure similar
to the Schmidt orthogonalization \{\cite{DH:CtoG}, pp. 27, 28 (3.1)-(3.4)\}.} 
of ${\cal G}(I)$:
\begin{equation}
a_J\equiv \vec{a}_{j_1}\wedge \vec{a}_{j_2}\wedge\ldots \wedge \vec{a}_{j_n},
\label{eq:P13.2}
\end{equation}
with $j_k=k$ or $0$ and elimination of elements with $j_k=0$.  Hence $j_1<j_2<\ldots <j_n$. 
$J$ stands for the combined 
index 
\begin{equation}
J\equiv (j_1,j_2, \ldots ,j_n).
\end{equation}

For all $j_k=0$ we define
\begin{equation}
a_J=a_{(0,0,\ldots,0)}\equiv 1.
\end{equation}
The $a^J$ are the corresponding reciprocal blades:
\begin{equation}
a^J\equiv \vec{a}^{j_n} \wedge\ldots\wedge \vec{a}^{j_2}\wedge \vec{a}^{j_1}.
\end{equation}
\more
The vectors with upper index are reciprocal vectors with
\begin{equation}
\vec{a}^k \vec{a}_m=\delta^k_m \,\,\,\,\mbox{(the Kronecker delta.)}
\end{equation}
For further details of the notation employed here compare~\cite{DH:CtoG}, pages 30 and 31. 
 The scalar differential operator $a_J \ast \partial_X$ after the 
last equality of (\ref{eq:P13.1a}) is not to be applied to $a^J,$ but rather to any multivector function to 
which the  multivector derivative operator $\partial_X$ on the left hand side is intended to be applied. 
\label{P13}
\epro
\prf{\ref{P13}}
See the definition of the algebraic properties of $\partial_X$ in definition \ref{D12}(i). For the last equality in (\ref{eq:P13.1a}) take into account that $a_J \ast \partial_X$ is 
algebraically scalar and~\cite{DH:CtoG}, chapter 1-1, (1.11).
\alles

\bpro[{\normalfont {\itshape constant scalar factor}}]
\normalfont
Another algebraic property of the multivector derivative $\partial_X$ is that we have for constant 
scalar factors $\lambda$:
\begin{align}
&\partial_X (\lambda F) = \lambda \,\,\partial_X F.
\end{align}
\label{P13a}
\epro
\prf{\ref{P13a}}
\begin{align}
&\partial_X(\lambda F)
\stackrel{\mathrm{(\ref{eq:P13.1a})}}{=}
\sum_J a^J a_J\ast\partial_X (\lambda F)
\stackrel{\mathrm{P. \ref{P5}}}{=}
\sum_J a^J \lambda \,\,a_J\ast\partial_X F 
\nonumber \\
&\stackrel{\mbox{\scriptsize \cite{DH:CtoG}, chap. 1-1, (1.11)}}{=}
\lambda\sum_J a^J  a_J\ast\partial_X F 
\stackrel{\mathrm{(\ref{eq:P13.1a})}}{=}
\lambda \,\,\partial_X F.
\end{align}

\bpro
\normalfont
\begin{align}
\partial_X &= \sum_{r=0}^n\partial_{\langle X \rangle_r} \mbox{ with } \nonumber \\
\partial_{\langle X \rangle_r} &=\langle \partial_X \rangle_r=
\sum_J \langle a^J \rangle_r \langle a_J \rangle_r\ast \partial_X.
\end{align}
$\partial_{\langle X \rangle_r}$ is thus the derivative with respect to a 
variable $\langle X \rangle_r \in {\cal G}^r(I).$
\label{P14}
\epro
\prf{\ref{P14}}
\begin{equation}
a_J\ast \partial_X \stackrel{\mbox{\scriptsize \cite{DH:CtoG}, p. 13 (1.46)}}{=} 
\sum_{r=0}^n \langle a_J \rangle_r\ast\partial_X
\stackrel{\mathrm{(ibidem)}}{=}
\sum_{r=0}^n \langle a_J \rangle_r\ast\langle \partial_X \rangle_r
\end{equation}
$\langle \partial_X \rangle_r$ is the $r$-blade part of $\partial_X,$ which has the 
algebraic properties of a multivector (Def. \ref{D12}, P. \ref{P13}). 
The sum 
$\sum_J a^J a_J\ast \partial_X$ is therefore naturally performed in two steps:
\begin{description}
\item[1.] Sum up over all index sets $J=(j_1,j_2,\ldots,j_n)$ with $r$ non-zero members.
\item[2.] Sum up over all $r=0 \ldots n$:
\end{description}
\begin{equation}
\partial_X =\sum_{r=0}^n \sum_J \langle a^J \rangle_r \langle a_J \rangle_r \ast \partial_X.
\end{equation}
\alles

\bpro
\begin{equation}
\partial=\sum_{r=0}^n\partial_{\bar{r}}.
\end{equation}
\label{P16}
\epro
\prf{\ref{P16}}
\begin{equation}
\partial_X \stackrel{\mbox{\scriptsize P. \ref{P14}}}{=} \sum_{r=0}^n\partial_{\langle X \rangle_r}
\stackrel{\mbox{\scriptsize P. \ref{P14}}}{=} \sum_{r=0}^n \langle \partial_{ X} \rangle_r
\stackrel{\mbox{\scriptsize Def. \ref{D15}}}{=}\sum_{r=0}^n\partial_{\bar{r}}.
\end{equation}

\bpro
\begin{equation}
\partial_X = \partial_A A\ast \partial_X.
\end{equation}
\label{P20}
\epro
\prf{\ref{P20}}
\begin{align}
\partial_A A\ast \partial_X &\stackrel{\mbox{\scriptsize P. \ref{P13}}}{=} 
\sum_J a^J a_J\ast\partial_A (A\ast \partial_X)
\nonumber \\
&\stackrel{\mbox{\scriptsize P. \ref{P19}}}{=}
\sum_J a^J P(a_J)\ast \partial_X=\sum_J a^J a_J\ast\partial_X
\stackrel{\mbox{\scriptsize P. \ref{P13}}}{=}\partial_X.
\end{align}

\bpro[{\normalfont {\itshape derivative from differential}}]
\normalfont
\begin{equation}
\partial F \equiv \partial_X F(X) =\partial_A\underline{F}(X,A)
=\partial_A F_A(X)\equiv \underline{\partial}\,\underline{F}.
\end{equation}
$\underline{\partial}$ means the derivative with respect to the differential argument $A$ of $\underline{F}.$
\label{P21}
\epro
\prf{\ref{P21}}
\begin{equation}
\partial_X F(X)\stackrel{\mbox{\scriptsize P. \ref{P20}}}{=}
\partial_A A\ast \partial_XF\stackrel{\mbox{\scriptsize Def. \ref{D1}}}{=}
\partial_A\underline{F}(X,A).
\end{equation}

\section{\cred Adjoint of a Multivector Function}

\bdf[{\normalfont {\itshape adjoint}}]
\normalfont
The adjoint of a multivector function $F$ is
\begin{equation}
\overline{F}=\overline{F}(A^{\prime})\equiv \underline{\partial}(\underline{F}\ast A^{\prime}),
\end{equation}
or explicitly
\begin{equation}
\overline{F}(A^{\prime})\equiv \partial_A \left[ \left\{ A\ast\partial_X F(X)\right\}\ast A^{\prime}\right]
\stackrel{\mbox{\scriptsize Def. \ref{D1}}}{=} \partial_A \left[ \underline{F}(X,A) \ast A^{\prime}\right].
\end{equation}
\label{D22}
\edf

\bpro
\normalfont
\begin{equation}
\overline{F}(A^{\prime})=\partial_X(F\ast A^{\prime}).
\end{equation}
or explicitly
\begin{equation}
\overline{F} (X,A^{\prime}) = \partial_X (F(X) \ast A^{\prime}).
\end{equation}
\label{P23}
\epro
\prf{\ref{P23}}
\begin{align}
\overline{F} (X,A^{\prime}) 
&\stackrel{\mbox{\scriptsize Def. \ref{D22}}}{=} 
\partial_A \left[ \left\{ (A\ast \partial_X)F(X)\right\}\ast A^{\prime}\right]
\stackrel
{\stackrel{\mbox{\scriptsize $A \ast \partial_X$}}{\mbox{\scriptsize alg. scalar}}}{=}
\partial_A(A\ast \partial_X)\left[{F}(X) \ast A^{\prime}\right]
\nonumber \\
&\stackrel{\mbox{\scriptsize P. \ref{P20}}}{=}  
\partial_X (F(X) \ast A^{\prime}).
\end{align}

\bpro[{\normalfont {\itshape common definition of adjoint}}]
\begin{equation}
B\ast \overline{F}(A)=\underline{F}(B)\ast A.
\end{equation}
\label{P24}
\epro
\prf{\ref{P24}}
\begin{equation}
B\ast \overline{F}(A)\stackrel{\mbox{\scriptsize P. \ref{P23}}}{=}
B \ast \partial_X F\ast A 
\stackrel
{\stackrel{\mbox{\scriptsize $B \ast \partial_X$}}{\mbox{\scriptsize algebraic scalar}}}{=}
(B \ast \partial_X F)\ast A
\stackrel{\mbox{\scriptsize Def. \ref{D1}}}{=}
\underline{F}(B)\ast A.
\end{equation}

\bpro
\begin{equation}
P(\overline{F}(A))=\overline{F}(A).
\end{equation}
\label{P25}
\epro
\prf{\ref{P25}}
\begin{equation}
P(\overline{F}(A))
\stackrel{\mbox{\scriptsize P. \ref{P23}}}{=}
P(\partial_X\underbrace{(F\ast A)}_{\mathrm{scalar}})
=P(\partial_X)(F\ast A)
\stackrel{\mbox{\scriptsize P. \ref{P13}}}{=}
\partial_X(F\ast A)
\stackrel{\mbox{\scriptsize P. \ref{P23}}}{=}
\overline{F}(A).
\end{equation}

\bpro
\begin{equation}
\langle\overline{F}(A) \rangle_r
= \underline{\partial}_{\,\bar{r}}(\underline{F}\ast A)=\partial_{\bar{r}}F\ast A.
\end{equation}
\label{P26}
\epro
\prf{\ref{P26}}
\begin{equation}
\langle\overline{F}(A) \rangle_r
\stackrel{\mbox{\scriptsize Def. \ref{D22}}}{=}
\langle \underline{\partial}\underbrace{\underline{F}\ast A}_{\mathrm{scalar}} \rangle_r=
\langle \partial_B \rangle_r \underline{F}(X,B)\ast A
\stackrel{\mbox{\scriptsize Def. \ref{D15}}}{=}
 \underline{\partial}_{\,\bar{r}}(\underline{F}\ast A),
\end{equation}
\begin{equation}
\langle\overline{F}(A) \rangle_r
\stackrel{\mbox{\scriptsize P. \ref{P23}}}{=}
\langle \partial_X(\underbrace{F(X)\ast A}_{\mathrm{scalar}}) \rangle_r=
\langle \partial_X \rangle_r ({F}(X)\ast A)
\stackrel{\mbox{\scriptsize Def. \ref{D15}}}{=}
\partial_{\bar{r}}F(X)\ast A.
\end{equation}

\bpro[{\normalfont {\itshape linearity of adjoint}}]
\normalfont
\begin{align}
\overline{F}(A+B)&=\overline{F}(A)+\overline{F}(B), \\
\overline{F}(\lambda A)&=\lambda \overline{F}(A), \,\,\,\,\,\, \mbox{if $\lambda=\langle \lambda \rangle$ scalar.}
\end{align}
\label{P27}
\epro
\prf{\ref{P27}}
\begin{align}
\overline{F}(A+B)
&\stackrel{\mbox{\scriptsize P. \ref{P23}}}{=}
\partial_X (F\ast (A+B)) 
\stackrel{\stackrel{\mbox{\scriptsize linearity of}}{\mbox{\scriptsize scalar product}}}{=}
\partial_X (F\ast A+F\ast B) \nonumber \\
&\stackrel{\mbox{\scriptsize P. \ref{P13}}}{=}
\sum_J a^J a_J\ast \partial_X (F\ast A+F\ast B) \nonumber \\
&\stackrel{\mbox{\scriptsize P. \ref{P4}, distributivity}}{=}
\sum_J a^J a_J\ast \partial_X (F\ast A)+\sum_J a^J a_J\ast \partial_X (F\ast B) \nonumber \\
&\stackrel{\mbox{\scriptsize P. \ref{P13}}}{=}
\partial_X(F\ast A)+\partial_X(F\ast B)
\stackrel{\mbox{\scriptsize P. \ref{P23}}}{=}
\overline{F}(A)+\overline{F}(B),
\end{align}
where I mean the distributivity of geometric multiplication with respect to addition as in~\cite{DH:CtoG}
p. 3, (1.4), (1.5).
\newline
\more
\begin{align}
\overline{F}(\lambda A)
&\stackrel{\mbox{\scriptsize P. \ref{P23}}}{=}
\partial_X F\ast (\lambda A)
\stackrel{\mbox{\scriptsize scalar product}}{=}
\partial_X(\lambda F\ast A) \nonumber \\
&\stackrel{\mbox{\scriptsize P. \ref{P13}}}{=}
\sum_J a^J a_J\ast \partial_X (\lambda F\ast A)
\stackrel{\mbox{\scriptsize P. \ref{P5}}}{=}
\sum_J a^J \lambda a_J\ast \partial_X ( F\ast A)  \nonumber \\
&\hspace{-0.6cm}\stackrel{\mbox{\scriptsize \cite{DH:CtoG} p. 4, (1.11)}}{=}
\lambda \sum_J a^J  a_J\ast \partial_X ( F\ast A) 
\stackrel{\mbox{\scriptsize P. \ref{P23}, P. \ref{P13}}}{=}
\lambda \overline{F}(A).
\end{align}
\alles

\section{\cred Differentiating Sums and Products, Changing Variables}

\bpro[{\normalfont {\itshape sum rule}}]
\begin{equation}
\partial_{\overline{r}}(F+G)=\partial_{\overline{r}}F+\partial_{\overline{r}}G.
\end{equation}
\label{P28}
\epro
\prf{\ref{P28}}
\begin{align}
\partial_{\overline{r}}(F+G)
&\stackrel{\mbox{\scriptsize P. \ref{P14}}}{=}
\sum_J \langle a^J\rangle_{\overline{r}}\langle a_J\rangle_{\overline{r}}\ast \partial_X(F+G) 
\nonumber \\
&\stackrel{\mbox{\scriptsize P. \ref{P6}}}{=}
\sum_J \langle a^J\rangle_{\overline{r}}\left\{ \langle a_J\rangle_{\overline{r}}\ast \partial_XF
+\langle a_J\rangle_{\overline{r}}\ast \partial_XG \right\}
\nonumber \\
&\hspace{-0.6cm}\stackrel{\mbox{\scriptsize distributivity}}{=}
\sum_J \langle a^J\rangle_{\overline{r}} \langle a_J\rangle_{\overline{r}}\ast \partial_XF
+ \sum_J \langle a^J\rangle_{\overline{r}} \langle a_J\rangle_{\overline{r}}\ast \partial_XG 
\nonumber \\
&\stackrel{\mbox{\scriptsize P. \ref{P14}}}{=}
\partial_{\overline{r}}F+\partial_{\overline{r}}G,
\end{align}
where I again mean the distributivity of geometric multiplication with respect to addition as in~\cite{DH:CtoG}
p. 3, (1.4), (1.5).
\newline
\alles

\bpro[{\normalfont {\itshape product rule}}]
\begin{equation}
\partial_{\overline{r}}(FG)=\dot{\partial}_{\overline{r}}\dot{F}G+\dot{\partial}_{\overline{r}}F\dot{G}.
\label{eq:P29.1}
\end{equation}
\label{P29}
\epro
\prf{\ref{P29}}
\begin{align}
\partial_{\overline{r}}(FG)
&\stackrel{\mbox{\scriptsize P. \ref{P14}}}{=}
\sum_J \langle a^J\rangle_{\overline{r}}\langle a_J\rangle_{\overline{r}}\ast \partial_X(FG)
\nonumber \\
&\stackrel{\mbox{\scriptsize P. \ref{P7}}}{=}
\sum_J \langle a^J\rangle_{\overline{r}}\big\{(\langle a_J\rangle_{\overline{r}}\ast \partial_XF)G+
F(\hspace{-0.1cm}\underbrace{\langle a_J\rangle_{\overline{r}}\ast \partial_X}_{\mbox{\scriptsize algebraic scalar}}\hspace{-0.1cm})G\big\}
\nonumber \\
&=
\sum_J \langle a^J\rangle_{\overline{r}}\left\{(\langle a_J\rangle_{\overline{r}}\ast \dot{\partial}_X)\dot{F}G+
(\langle a_J\rangle_{\overline{r}}\ast \dot{\partial}_X)F\dot{G}\right\}
\nonumber \\
&\hspace{-0.6cm}\stackrel{\mbox{\scriptsize distributivity}}{=}
\sum_J \langle a^J\rangle_{\overline{r}}(\langle a_J\rangle_{\overline{r}}\ast \dot{\partial}_X)\dot{F}G+
\sum_J \langle a^J\rangle_{\overline{r}}(\langle a_J\rangle_{\overline{r}}\ast \dot{\partial}_X)F\dot{G}
\nonumber \\
&\stackrel{\mbox{\scriptsize P. \ref{P14}}}{=}
\dot{\partial}_{\overline{r}}\dot{F}G+\dot{\partial}_{\overline{r}}F\dot{G},
\end{align}
where I again mean the distributivity of geometric multiplication with respect to addition as in~\cite{DH:CtoG}
p. 3, (1.4), (1.5).
\newline
\alles

\bdf[{\normalfont {\itshape product rule variation}}]
\normalfont
Left \textit{and} right side derivation is indicated by:
\begin{equation}
\dot{F}\dot{\partial}_{\overline{r}}\, \dot{G}=\dot{F}\dot{\partial}_{\overline{r}}\, G
+F \dot{\partial}_{\overline{r}}\, \dot{G}.
\label{eq:D30.1}
\end{equation}
\label{D30}
\edf

\brem
\normalfont
Expanding the variation (\ref{eq:D30.1}) of the product rule explicitely in the $a_J$ blade basis (\ref{eq:P13.2})
of the geometric algebra ${\cal G}$ shows how (\ref{eq:P29.1}) and (\ref{eq:D30.1}) are algebraically different
(compare Proof \ref{P29}):
\begin{align}
\dot{F}\dot{\partial}_{\overline{r}}\, \dot{G}
&\stackrel{\mbox{\scriptsize D. \ref{D30}}}{=}
\dot{F}\dot{\partial}_{\overline{r}}\, G
+F \dot{\partial}_{\overline{r}}\, \dot{G}
\nonumber \\
&\stackrel{\mbox{\scriptsize P. \ref{P14}}}{=}
\sum_J  (\langle a_J\rangle_{\overline{r}}\ast \dot{\partial}_X)\dot{F}\langle a^J\rangle_{\overline{r}} G+
\sum_J  F\langle a^J\rangle_{\overline{r}} (\langle a_J\rangle_{\overline{r}}\ast \dot{\partial}_X)\dot{G}
\end{align}
\label{R31}
\erem

\bpro[{\normalfont {\itshape change of variables}}]
\normalfont
For $F(X)=G(f(X)),$ i.e. $f:X\rightarrow X^{\prime} =f(X)$
\begin{align}
\partial_{\overline{r}}\,F
 &= \dot{\partial}_{\overline{r}}\, G(\dot{f})=\langle\overline{f}(\dot{\partial}^{\prime})  \rangle_r \dot{G},
\\
\mbox{i.e. } \,\,\,\,\,\,\,\,\,
\partial_{\overline{r}}
&=\langle\partial_X  \rangle_r=
\langle\overline{f}({\partial}_{X^{\prime}})  \rangle_r=
\langle\overline{f}({\partial}^{\prime})  \rangle_r.
\end{align}
\label{P31}
\epro
\prf{\ref{P31}}
\begin{align}
\partial_{\overline{r}}\,F(X)
&\stackrel{\mbox{\scriptsize P. \ref{P21}}}{=}
\partial_A (A \ast \langle \partial_X  \rangle_r)G(f(X))
\nonumber \\
&\hspace{-0.6cm}
\stackrel{\mbox{\scriptsize \cite{DH:CtoG} p. 13, (1.46)}}{=}
\partial_{A_{\overline{r}}}\,(A_{\overline{r}}\ast \langle\partial_X  \rangle_r)G(f(X))
\nonumber \\
&\stackrel{\mbox{\scriptsize Def. \ref{D1}}}{=}
\partial_{A_{\overline{r}}}\,\partial_{\tau}\left. G(f(X+\tau P(A_{\overline{r}})))\right|_{\tau=0}
\nonumber \\
&\stackrel{\mbox{\scriptsize P. \ref{P8}}}{=}
\partial_{A_{\overline{r}}}\,\partial_{\tau}
\left. G(f(X)+\tau (A_{\overline{r}}\ast \langle\partial_X  \rangle_r)f(X))\right|_{\tau=0}
\nonumber \\
&\stackrel{\mbox{\scriptsize Def. \ref{D1}}}{=}
\partial_{A_{\overline{r}}}\,\left(\left\{A_{\overline{r}}\ast \langle\partial_X  \rangle_rf(X)  \right\}
\ast {\partial}_{X^{\prime}}\right)\left. G(X^{\prime})\right|_{X^{\prime}=f(X)}
\nonumber \\
&\stackrel{\mbox{\scriptsize Def. \ref{D1}}}{=}
\partial_{A_{\overline{r}}} 
\left( \underline{f}(A_{\overline{r}}) \ast {\partial}_{X^{\prime}}
\right)\left. G(X^{\prime})\right|_{X^{\prime}=f(X)}
\nonumber \\
&\stackrel{\mbox{\scriptsize P. \ref{P26}}}{=}
\langle\overline{f}({\partial}_{X^{\prime}})  \rangle_r \left. G(X^{\prime})\right|_{X^{\prime}=f(X)}.
\end{align}
\newline
\alles

\brem
\normalfont
Remark {\ref{R10}} also applies to proof {\ref{P31}}. I conclude from proof {\ref{P31}} that formulas
(2.26a+b) in~\cite{DH:CtoG}, p. 56 are slightly incorrect. There the grade selector $\langle \ldots  \rangle_r$ 
is applied to the argument of $\overline{f}$, but in order to be correct it should be applied to $\overline{f}$
itself, i.e. appear around  $\overline{f}$ as in P. \ref{P31}.  
\erem

\bpro[{\normalfont {\itshape using the full derivative $\partial_X$}}]
\normalfont
\begin{align}
&\hspace{-1.5cm}\mbox{Sum rule:} &\partial(F+G)=\partial F+ \partial G, \\
&\hspace{-1.5cm}\mbox{Constant scalar factor $\lambda$:} &\partial(\lambda F)=\lambda \partial F, 
\label{eq:P32.2}\\
&\hspace{-1.5cm}\mbox{Product rule:} &\partial(FG)=\dot{\partial}\dot{F}G+\dot{\partial}F\dot{G}, \\
&\hspace{-1cm}\mbox{with variation:} &\dot{F}\dot{\partial}\dot{G}=\dot{F}\dot{\partial}{G}+{F}\dot{\partial}\dot{G}.
\end{align}
For $F(X)=G(f(X)),$ i.e. $f: X\rightarrow X^{\prime}=f(X)$
\begin{align}
\mbox{Chain rule: }\hspace{3.4cm}\partial F &= \dot{\partial}G(\dot{f}) = \overline{f}(\dot{\partial}) \dot{G}, \\
\mbox{i.e.} \hspace{3.4cm}\,\,\,\,\,\,  \partial &= \overline{f}(\partial_{X^{\prime}})=\overline{f}({\partial^{\prime}}).
\end{align}
\label{P32}
\epro
\prf{\ref{P32}}
One just needs to take the sum over all grades $r:$ $\sum_{r=0}^n$ on both left and right hand sides of P. \ref{P28},
P. \ref{P29}, Def. \ref{D30} and P. \ref{P31}, as well as take into account P. \ref{P14} and that
\begin{equation}
\overline{f}(\partial^{\prime})=\sum_{r=0}^n\langle\overline{f}(\partial^{\prime})   \rangle_r, \hspace{1cm}
 \mbox{compare~\cite{EH:VDcalc}, (13).}
\end{equation}
For (\ref{eq:P32.2}) compare P. \ref{P13a}.

\brem
\normalfont
The above results obtain, if $F=F(X)$ is defined for $X=X_{\overline{r}}=\langle X \rangle_r\in {\cal G}^r(I).$
The results can be adapted to functions on any linear subspace of ${\cal G}(I)$. The differential 
$\underline{F}=\underline{F}(X,A)$ (Def. \ref{D1}) and the adjoint $\overline{F}=\overline{F}(X,A^{\prime})$
(Def. \ref{D22}) are the only two linear functions of $A$ (and $A^{\prime}$ respectively) that can be formed from 
$F$ using the derivative $\partial=\partial_X$ and the scalar product (compare the definitions Def. \ref{D1} 
and Def. \ref{D22}.)
\label{R33}
\erem

\section{\cred The Scalar Case}

\bpro[{\normalfont {\itshape scalar differential calculus}}]
\normalfont
For a multivector function of a scalar variable $X=X(\tau)=f(\tau):$
\begin{align}
&\label{P34.1}\mbox{Multivector derivative}  \hspace{0.7cm} \partial_{\tau}X=\langle \partial_{\tau}\rangle X=\frac{dX}{d\tau},   \\
&\label{P34.2}\mbox{Differential}\hspace{2.3cm}   \underline{X}(\lambda)=\lambda\ast \partial X 
= \lambda \partial_{\tau} X=\lambda \frac{dX}{d\tau}, \\
&\label{P34.3}\mbox{Adjoint} \hspace{2.8cm}   \overline{X}(A)=\partial_{\tau}X\ast A = \left(\frac{dX}{d\tau}\right)\ast A,  \\
&\label{P34.4}\mbox{Chain rule}  \hspace{2.3cm}   \frac{dF}{d\tau}= \overline{f}(\partial_{X})F(X)=\frac{dX}{d\tau}\ast \partial_X F(X).
\end{align}
\more
The second expression of equation (2.27d) in \cite{DH:CtoG} seems to be sligthly wrong compared to (\ref{P34.4}).
In the special case of a scalar function $X=x(\tau)=\langle x(\tau)\rangle,$ only the scalar part 
$\alpha=\langle A\rangle$ of $A$ contributes to the adjoint:
\begin{equation}
\overline{X}(\alpha)=\alpha \frac{dx}{d\tau},
\label{P34.5}
\end{equation}
and the chain rule for such a scalar function has the form:
\begin{equation}
\frac{dF}{d\tau}=\frac{dx}{d\tau} \frac{dF}{dx}.
\end{equation}
\label{P34}
\epro
\more
\prf{\ref{P34}}\\
For the multivector derivative (\ref{P34.1}):\\
 Because
$\tau$ is scalar, $\langle \partial_{\tau}\rangle_r X=0$ for $r\neq 0,$
$\Rightarrow \partial=\partial_{\tau}=\frac{d}{d\tau}$ \\
$\Rightarrow \mbox{for the differential (\ref{P34.2})}:$
\begin{equation}
\underline{X}(A)
\stackrel{\mbox{\scriptsize Def. \ref{D1}}}{=}
A\ast \partial X = (A\ast \partial_{\tau})X
\stackrel{\mbox{\scriptsize scalar product}}{=}
\langle A \rangle \partial_{\tau} X
\stackrel{\mbox{\scriptsize $\lambda \equiv\langle A \rangle $}}{=}
\lambda \partial_{\tau} X.
\end{equation}
For the adjoint (\ref{P34.3}) compare P. {\ref{P23}}. \\
For the chain rule (\ref{P34.4}) compare P. \ref{P31} and (\ref{P34.3}) with $A=\partial_X$.\\
We further have from (\ref{P34.3}) for scalar $x(\tau)=f(\tau):$ $\frac{dx}{d\tau}=\langle \frac{dx}{d\tau}\rangle$
\begin{equation}
\overline{x}(A) = \frac{dx}{d\tau}\ast A 
\stackrel{\mbox{\scriptsize scalar product}}{=}
\frac{dx}{d\tau}\ast \langle A\rangle
\stackrel{\mbox{\scriptsize $\alpha \equiv\langle A \rangle $}}{=}
\frac{dx}{d\tau}\alpha
\end{equation}
and from (\ref{P34.4})
\begin{equation}
\Rightarrow \frac{dF}{d\tau}=\overline{f}(\partial_{x})F(x)=
\langle \frac{dx}{d\tau}\rangle\ast \partial_x F(x)=\frac{dx}{d\tau}\frac{d}{dx}F(x).
\end{equation}
\alles

\brem
\normalfont
In scalar differential calculus differentials (\ref{P34.2}) and adjoints (\ref{P34.5}) are identical
for $\lambda = \alpha$. 
The distinction with the derivative (\ref{P34.1})
is trivial (i.e. only the scalar factor $\lambda$). 
The single concept of derivative in elementary differential calculus is therefore
 now generalized to three distinct, 
related concepts of multivector derivative, differential and adjoint.
\erem

\section{\cred Basic Multivector Derivatives \label{sc:BMD}}

In the following I assume
\begin{description}
\item[(i)] $X=F(X)$ to be the identity function on some linear subspace of ${\cal G}(I)$ of dimension $d,$
\item[(ii)] $\Slash{A} \equiv P_{\mathrm{subpace}}(A)$ to be the projection into the above mentioned 
$d-$di\-men\-sion\-al
subspace of ${\cal G}(I),$ 
\item[(iii)] that singularities at $X=0$ are to be excluded.
\end{description}

\bpro
\begin{equation}
A\ast\partial_X X = \dot{\partial}_X\dot{X}\ast A =\Slash{A}.
\label{eq:P36.1a}
\end{equation}
\label{P36}
\epro

\brem
\normalfont
To help avoid confusion I refer to Remarks \ref{R2} and \ref{R13} in order to clarify that the 
computation of (\ref{eq:P36.1a}) implies the following brackets:
\begin{align}
(A\ast\partial_X) X = \dot{\partial}_X (\dot{X}\ast A) =\Slash{A}.
\end{align}
\erem

\prf{\ref{P36}}
\begin{align}
&F(X)\equiv X, \nonumber \\
&A\ast\partial_X X 
\stackrel{\mbox{\scriptsize Def. \ref{D1}}}{=}
\underline{F}(A) 
\stackrel{\mbox{\scriptsize P. \ref{P2}}}{=}
\underline{F}(P(A))=\underline{F}(\Slash{A})
\stackrel{\mbox{\scriptsize P. \ref{P18}}}{=}\Slash{A},
\label{P36.1} \\
&\dot{\partial}_X\dot{X}\ast A 
\stackrel{\mbox{\scriptsize Def. \ref{D12}}}{=}
\sum_{J_d}a^J \underbrace{a_J \ast \dot{\partial}_X\dot{X}}_{\mbox{\scriptsize using (\ref{P36.1})}}\ast A
= \sum_{J_d}a^J a_J\ast A=\Slash{A},
\end{align}
where $J_d$ is the index subset for the multivector base of the assumed linear 
$d-$di\-men\-sion\-al subspace of ${\cal G}(I)$ of $X.$

\bpro
\normalfont
\begin{equation}
A\ast\partial_X \tilde{X} = \dot{\partial}_X\dot{\tilde{X}}\ast A =\Slash{\tilde{A}},
\end{equation}
where the tilde operation $(\tilde{\,\,\,})$ indicates \textit{reversion} 
as defined in~\cite{DH:CtoG}, pp. 5,6, (1.17) to (1.20).
\label{P37}
\epro

\prf{\ref{P37}}
\begin{align}
&\hspace{-0.4cm}\underbrace{A\ast\partial_X}_{\mbox{\scriptsize algebraic scalar}} \hspace{-0.4cm}\tilde{X} 
=\widetilde{(A\ast\partial_X X)}
\stackrel{\mbox{\scriptsize P. \ref{P36}}}{=} \Slash{\tilde{A}}  \\
&\hspace{0.2cm}\Slash{\tilde{A}}
\stackrel{\mbox{\scriptsize P. \ref{P36}}}{=}
\dot{\partial}_X\dot{X}\ast \tilde{A}\,\,\,
\stackrel{\mbox{\scriptsize \cite{DH:CtoG} p. 13, (1.48)}}{=}
\dot{\partial}_X\dot{\tilde{X}}\ast A
\end{align}

\bpro
\begin{equation}
\partial_XX=d.
\end{equation}
\label{P38}
\epro
\prf{\ref{P38}}
\begin{equation}
\partial_XX\stackrel{\mbox{\scriptsize Def. \ref{D12}}}{=}
\sum_{J_d} a^J a_J\ast\partial_XX
\stackrel{\mbox{\scriptsize P. \ref{P36}}}{=}
\sum_{J_d} a^J a_J=\sum_{J_d} \delta^{J}{}_{J}=d,
\end{equation}
where $J_d$ is the index subset for the multivector base of the assumed linear 
$d-$di\-men\-sion\-al subspace of ${\cal G}(I)$ of $X.$

\bpro
\begin{align}
\partial_X|X|^2= 2 \tilde{X}.
\end{align}
\label{P39}
\epro
\prf{\ref{P39}}
\begin{align}
\partial_X|X|^2
&\stackrel{\mbox{\scriptsize \cite{DH:CtoG} p. 13, (1.49)}}{=}
\partial_X \langle X\tilde{X}  \rangle
\stackrel{\mbox{\scriptsize Def. \ref{D12}, P. \ref{P13}}}{=}
\sum_{J_d} a^J a_J\ast\partial_X\langle X\tilde{X}  \rangle \nonumber
\\
&\stackrel{\stackrel{\mbox{\scriptsize P. \ref{P7}, P. \ref{P3}}}{\mbox{\scriptsize distributivity}}}{=}
\sum_{J_d} a^J a_J\ast\dot{\partial}_X\langle \dot{X}\tilde{X}  \rangle+
\sum_{J_d} a^J a_J\ast\dot{\partial}_X\langle X\dot{\tilde{X}}  \rangle
 \\
&\hspace{0.4cm}\stackrel{\mbox{\scriptsize Def. \ref{D12}, P. \ref{P13}}}{=}
\dot{\partial}_X\langle \dot{X}\tilde{X}  \rangle+
\dot{\partial}_X\langle X\dot{\tilde{X}}  \rangle
=
\dot{\partial}_X( \dot{X}\ast\tilde{X}  )+
\dot{\partial}_X(X\ast\dot{\tilde{X}}  )
\nonumber \\
&\stackrel{\mbox{\scriptsize \cite{DH:CtoG} p. 13, (1.47a)}}{=}
\dot{\partial}_X( \dot{X}\ast\tilde{X}  )+
\dot{\partial}_X(\dot{\tilde{X}} \ast X )
\stackrel{\mbox{\scriptsize P. \ref{P36},P. \ref{P37}}}{=}
\tilde{X}+\tilde{X}=2\tilde{X}, \nonumber
\end{align}
where I mean the distributivity of geometric multiplication with respect to addition as in~\cite{DH:CtoG}
p. 3, (1.4), (1.5). $J_d$ is the index subset for the multivector base of the assumed linear 
$d-$di\-men\-sion\-al subspace of ${\cal G}(I)$ of $X.$

\bpro
\begin{equation}
A\ast \partial_X X^k=\Slash{A} X^{k-1}+X\Slash{A} X^{k-2}+\ldots +X^{k-1}\Slash{A}.
\end{equation}
\label{P40}
\epro
\prf{\ref{P40}}
\begin{align}
&\,\,\,F \equiv X, \nonumber \\
A\ast \partial_X &X^k =A\ast \partial_X F^k
\stackrel{\mbox{\scriptsize Def. \ref{D1}}}{=}
\underline{(F^k)}
\stackrel{\mbox{\scriptsize P. \ref{P7}}}{=}
\underline{F}F^{k-1}+F \underline{F} F^{k-2}+\ldots + F^{k-1}\underline{F} \nonumber \\
&\!\!\!\stackrel{\mbox{\scriptsize P. \ref{P36}, $F\equiv X$}}{=}
\Slash{A} X^{k-1}+X\Slash{A} X^{k-2}+\ldots +X^{k-1}\Slash{A}.
\end{align}

\bpro
\begin{equation}
\partial_X |X|^k=k|X|^{k-2}\tilde{X}.
\end{equation}
\label{P41}
\epro
\prf{\ref{P41}}
\begin{align}
&f(X)\equiv |X|^2, \nonumber \\
&\partial_X |X|^k=\partial_X (|X|^2)^{\frac{k}{2}}
\stackrel{\mbox{\scriptsize P. \ref{P31}}}{=}
\overline{f}(\partial_{\alpha})\left. \alpha^{\frac{k}{2}}\right|_{\alpha=f(X)=|X|^2} \nonumber \\
&\hspace{1cm}\stackrel{\mbox{\scriptsize P. \ref{P23}}}{=}
\dot{\partial}_X(|\dot{X}|^2\ast \partial_{\alpha})\left. \alpha^{\frac{k}{2}}\right|_{\alpha=f(X)=|X|^2}
=(\partial_X|X|^2)\partial_{\alpha}\left. \alpha^{\frac{k}{2}}\right|_{\alpha=f(X)=|X|^2}\nonumber \\
&\hspace{1cm}
\stackrel{\mbox{\scriptsize P. \ref{P39}}}{=}
2\tilde{X}\frac{k}{2}|X|^{k-2}=
k|X|^{k-2}\tilde{X}.
\end{align}

\bpro
\begin{equation}
\partial_X \log |X| = \frac{\tilde{X}}{|X|^2}.
\end{equation}
\label{P42}
\epro
\prf{\ref{P42}}
\begin{align}
&f(X)\equiv |X|, \nonumber \\
&\partial_X \log |X|
\stackrel{\mbox{\scriptsize P. \ref{P31}}}{=}
\overline{f}(\partial_{\alpha}) \left. \log \alpha\right|_{\alpha=f(X)=|X|}
\stackrel{\mbox{\scriptsize P. \ref{P23}}}{=}
\dot{\partial_X}(|\dot{X}|\ast\partial_{\alpha})\left. \log \alpha\right|_{\alpha=f(X)=|X|}
\nonumber \\
&\hspace{1.7cm}= (\partial_X |X|)\frac{1}{|X|}
\stackrel{\mbox{\scriptsize P. \ref{P41}}}{=}
|X|^{-1}\tilde{X} \frac{1}{|X|}= \frac{\tilde{X}}{|X|^2}.
\end{align}

\bpro
\begin{equation}
A\ast \partial_X(|X|^k X)=|X|^k \left(\Slash{A}+k\frac{A\ast \tilde{X}X}{|X|^2}\right).
\end{equation}
\label{P43}
\epro
\prf{\ref{P43}}
\begin{align}
&A\ast \partial_X(|X|^k X)
\stackrel{\mbox{\scriptsize P. \ref{P7}}}{=}
(A\ast \partial_X |X|^k)X+ |X|^k A\ast \partial_X X
 \\
&\stackrel{\mbox{\scriptsize Def. \ref{D12}(ii), P. \ref{P41}, P. \ref{P36} }}{=}
A\ast (\underbrace{k |X|^{k-2}}_{\mathrm{scalar}}\tilde{X})X + |X|^k\Slash{A}
=|X|^k \left(\Slash{A}+k\frac{A\ast \tilde{X}X}{|X|^2}\right).
\nonumber
\end{align}

\bpro
\begin{equation}
\partial_X(|X|^kX) = |X|^k \left(d+k\frac{\tilde{X}X}{|X|^2}   \right).
\end{equation}
\label{P44}
\epro
\prf{\ref{P44}}
\begin{align}
\hspace{-1cm}\partial_X(|X|^kX)
&\stackrel{\mbox{\scriptsize P. \ref{P29}, Rem. {\ref{R33}}}}{=}
(\partial_X |X|^k)X +\dot{\partial}_X\underbrace{|X|^k}_{\mathrm{scalar}} \dot{X}
\nonumber \\
&\hspace{0.6cm}
\stackrel{\mbox{\scriptsize P. \ref{P41}}}{=}
(k |X|^{k-2}\tilde{X})X +|X|^k \partial_X X
\nonumber \\
&\hspace{0.6cm}\stackrel{\mbox{\scriptsize P. \ref{P38}}}{=}
k|X|^k \frac{\tilde{X}X}{|X|^2} + |X|^k d=
|X|^k \left(d+k\frac{\tilde{X}X}{|X|^2}   \right).
\end{align}

\bpro
\normalfont
For $X=\sum_{r=0}^{n}X_{\bar{r}}$ defined on the whole of ${\cal G}(I)$ and $\partial=\partial_X$
\begin{equation}
A\ast \partial X = \dot{\partial} \dot{X}\ast A = P(A).
\end{equation}
\label{P45}
\epro
\prf{\ref{P45}}
In Proposition \ref{P36}, the assumed "subspace" becomes now all of ${\cal G}(I)$, therefore
\Slash{A}=P(A).

\bpro
\normalfont
For $X=\sum_{r=0}^{n}X_{\bar{r}}$ defined on the whole of ${\cal G}(I)$ and $\partial=\partial_X$
\begin{equation}
A \ast \partial_{\bar{r}}X=\dot{\partial}_{\bar{r}} \dot{X}\ast A =P(A_{\bar{r}}).
\end{equation}
\label{P46}
\epro
\prf{\ref{P46}}
\begin{equation}
A \ast \partial_{\bar{r}}X
\stackrel{\mbox{\scriptsize \cite{DH:CtoG}, p. 13, (1.45a)}}{=}
A_{\bar{r}}\ast \partial_X X \stackrel{\mbox{\scriptsize P. \ref{P45}}}{=}
P(A_{\bar{r}}),
\label{eq:P46.1}
\end{equation}
and
\begin{align}
\dot{\partial}_{\bar{r}} \dot{X}\ast A
&\stackrel{\mbox{\scriptsize P. \ref{P14}}}{=}
\sum_{J}\langle a^J \rangle_r\underbrace{\langle a_J \rangle_r\ast\dot{\partial}_X \dot{X}}\ast A
\stackrel{\mbox{\scriptsize (\ref{eq:P46.1})}}{=}
\sum_{J}\langle a^J \rangle_r P(\langle a_J \rangle_r)\ast A
\nonumber \\
&\hspace{-0.8cm}\stackrel{\mbox{\scriptsize \cite{DH:CtoG}, p. 13, (1.45a)}}{=}
\sum_{J}\langle a^J \rangle_r \langle a_J \rangle_r\ast A_{\bar{r}}=P(A_{\bar{r}}).
\end{align}

\bpro
\begin{equation}
\partial_{\bar{r}}X=\partial X_{\bar{r}} = \partial_{\bar{r}} X_{\bar{r}}=
\left(
\begin{array}{c}
n \\
r
\end{array}
\right).
\end{equation}
\label{P47}
\epro
\prf{\ref{P47}}
\begin{align}
\partial_{\bar{r}}X
&\stackrel{\mbox{\scriptsize P. \ref{P14}}}{=}
\sum_{J}\langle a^J \rangle_r\langle a_J \rangle_r\ast\partial_X X  
\stackrel{\mbox{\scriptsize \cite{DH:CtoG}, p. 13, (1.45a)}}{=}
\sum_{J}\langle a^J \rangle_r\underbrace{\langle a_J \rangle_r\ast\partial_{\langle X\rangle_r} X}
\nonumber \\
&\stackrel{\mbox{\scriptsize P. \ref{P46}}}{=}
\sum_{J}\langle a^J \rangle_r P(\langle a_J \rangle_r)=\sum_{J}\langle a^J \rangle_r \langle a_J \rangle_r
=\left(
\begin{array}{c}
n \\
r
\end{array}
\right),
\end{align}
where 
$\left(
\begin{array}{c}
n \\
r
\end{array}
\right)$
is the dimension of the $r$-vector subspace of ${\cal G}(I)$. 
\newline
\more
${\partial}_{\bar{r}}{X}_{\bar{r}}$ is in this
subspace and we therefore have according to P. \ref{P14} that
\begin{align}
{\partial}_{\bar{r}}{X}_{\bar{r}}
\stackrel{\mbox{\scriptsize P. \ref{P14}}}{=}
\sum_{J}\langle a^J \rangle_r \langle a_J \rangle_r \ast {\partial}_{\bar{r}}{X}_{\bar{r}}
\stackrel{\mbox{\scriptsize P. \ref{P36}}}{=}
\sum_{J}\langle a^J \rangle_r \langle a_J \rangle_r
=
\left(
\begin{array}{c}
n \\
r
\end{array}
\right)\,\, ,
\end{align}
and
\begin{align}
\partial X_{\bar{r}}
&\stackrel{\mbox{\scriptsize P. \ref{P14}}}{=}
\sum_s \partial_{\bar{s}} X_{\bar{r}} \stackrel{\mbox{\scriptsize P. \ref{P14}}}{=}
\sum_s \sum_J \langle a^J \rangle_s \langle a_J \rangle_s \ast \partial_{\bar{s}} X_{\bar{r}} 
\nonumber \\
&\stackrel{\mbox{\scriptsize P. \ref{P3}}}{=}
\sum_s \sum_J \langle a^J \rangle_s \left\langle \langle a_J \rangle_s \ast \partial_{\bar{s}} X \right\rangle_r
\stackrel{\mbox{\scriptsize P. \ref{P46}}}{=}
\sum_s \sum_J \langle a^J \rangle_s \left\langle P(\langle a_J \rangle_s) \right\rangle_r
\nonumber \\
&=\sum_s \sum_J \langle a^J \rangle_s \underbrace{\left\langle \langle a_J \rangle_s \right\rangle_r}_{
=\delta_{sr}\langle a_J \rangle_s}
=\sum_J \langle a^J \rangle_r \langle a_J \rangle_r 
=\left(
\begin{array}{c}
n \\
r
\end{array}
\right).
\end{align}
\alles

\bpro
\beqn
\partial X = 
\sum_r \left(
\begin{array}{c}
n \\
r
\end{array}
\right)
= 2^n.
\end{equation}
\label{P48}
\epro
\prf{\ref{P48}}
\beqn
\partial X
\stackrel{\mbox{\scriptsize P. \ref{P14}}}{=}
\sum_r \partial_{\bar{r}} X 
\stackrel{\mbox{\scriptsize P. \ref{P47}}}{=}
\sum_r 
\left(
\begin{array}{c}
n \\
r
\end{array}
\right)
= 2^n.
\end{equation}

\bpro
\normalfont
For $A=P(A)$, $K=\frac{1}{2}(r+s-|r-s|)$ and 
$\left(
\begin{array}{c}
i \\
j
\end{array}
\right) =0$ if $j>i$
\beqn
\partial_{\bar{s}}\langle X A_{\bar{r}}\rangle_m
=
\langle A_{\bar{r}}\partial_{\bar{s}} \rangle_m X
=
\left(
\begin{array}{c}
r \\
k
\end{array}
\right)
\left(
\begin{array}{c}
n-r \\
s-k
\end{array}
\right)
\delta^m_{r+s-2k}
A_{\bar{r}},
\end{equation}
with $\delta^m_{r+s-2k}$ the Kronecker symbol for non-negative integers.
\label{P49}
\epro
\prf{\ref{P49}}
(This step by step proof extends over the next two pages.)
I will now assume $A_r$ to be a simple $r$-blade $A_r=\vec{a}_1 \vec{a}_2 \ldots \vec{a}_r,$ i.e. a geometric product
of $r$ orthogonal unit vectors. We now select an orthonormal basis of ${\cal G}(I)$ which 
includes $\vec{a}_1, \vec{a}_2,\ldots, \vec{a}_r.$ The first step in our proof is to show the
formula
\begin{equation}
\partial_{\bar{s}}\langle X A_r\rangle_m = \partial_{\bar{s}}\langle X_{\bar{s}} A_r\rangle_m.
\label{eq:P49.B}
\end{equation}
\newline
\more
The left hand side of (\ref{eq:P49.B}) yields
\begin{align}
\partial_{\bar{s}}\langle X A_r\rangle_m
&\stackrel{\mbox{\scriptsize P. \ref{P14}}}{=}
\sum_J \langle a^J \rangle_s \left\langle 
\underbrace{\langle a_J \rangle_s \ast \partial_{\bar{s}} X} A_r \right\rangle_m
\nonumber \\
&\stackrel{\mbox{\scriptsize P. \ref{P46}}}{=}
\sum_J \langle a^J \rangle_s \left\langle \langle a_J \rangle_s A_r \right\rangle_m,
\label{eq:P49.BL}
\end{align}
and the right hand side of  (\ref{eq:P49.B}) gives
\begin{align}
\partial_{\bar{s}} \langle X_{\bar{s}} A_r\rangle_m
&\stackrel{\mbox{\scriptsize P. \ref{P14}}}{=}
\sum_J \langle a^J \rangle_s \left\langle 
\langle a_J \rangle_s \ast \underbrace{ \partial_{\bar{s}} X_{\bar{s}} } A_r \right\rangle_m
\nonumber \\
&\stackrel{\mbox{\scriptsize P. \ref{P46}}}{=}
\sum_J \langle a^J \rangle_s 
\left\langle \dot{\partial}_{\bar{s}} \dot{X} \ast \langle a_J \rangle_s A_r \right\rangle_m
\nonumber \\
&\stackrel{\mbox{\scriptsize P. \ref{P46}}}{=}
\sum_J \langle a^J \rangle_s \left\langle
\langle a_J \rangle_s \ast \partial_{\bar{s}} X A_r \right\rangle_m
\stackrel{\mbox{\scriptsize (\ref{eq:P49.BL})}}{=}
\partial_{\bar{s}} \langle X A_r\rangle_m,
\end{align}
where to obtain the second equality we first set in P. \ref{P46} 
$r=s$ and $A=\langle a_J \rangle_s:$
\begin{align}
\langle a_J \rangle_s\ast \partial_{\bar{s}}X = \dot{\partial}_{\bar{s}}\dot{X}\ast \langle a_J \rangle_s.
\label{eq:P49.extra}
\end{align}
Second, by applying the $s$-grade selector $\langle \,\,\, \rangle_s$ on both sides of (\ref{eq:P49.extra}) 
we get the necessary relationship
\begin{align}
\langle a_J \rangle_s\ast \partial_{\bar{s}}X_{\bar{s}} 
= \dot{\partial}_{\bar{s}}
\underbrace{\dot{X}\ast \langle a_J \rangle_s}_{\mathrm{alg. scalar}}.
\end{align}
This completes the proof of (\ref{eq:P49.B}). 
Writing the vector factors $\langle a^J \rangle_s,$  $\langle a_J \rangle_s$ and $A_r$ in the last
expression of (\ref{eq:P49.BL}) explicitely, we obtain from (\ref{eq:P49.B}) that
\begin{align}
\partial_{\bar{s}}\langle X A_r\rangle_m 
&= \partial_{\bar{s}}\langle X_{\bar{s}} A_r\rangle_m
= \sum_J \langle a^J \rangle_s \left\langle \langle a_J \rangle_s A_r \right\rangle_m
\nonumber \\
&= \sum_{j_1<\ldots <j_s} \vec{a}^{j_{s}}\ldots\vec{a}^{j_2}\vec{a}^{j_1}
\langle 
 \underbrace{
  \underbrace{\vec{a}_{j_1}\vec{a}_{j_2}\ldots\vec{a}_{j_s}}_{\mbox{\scriptsize $s$ vectors}} 
  \underbrace{\vec{a}_1\vec{a}_2\ldots\vec{a}_r}_{\mbox{\scriptsize $r$ vectors}} 
 }_{\mbox{\scriptsize $r+s$ vectors}}
\rangle_m,
\end{align}
where all vectors $\vec{a}_{j_1}, \vec{a}_{j_2}, \ldots,\vec{a}_{j_s}$ are taken from the aforementioned 
orthonormal basis of ${\cal G}(I)$ and $j_s \leq n$.
\newline
\more
Selecting the $m$-grade part (i.e. the $m=r+s-2k$-grade part) from a product of $r+s$ vectors means, that because
of the orthogonality (all $s$ vectors $\vec{a}_{j_1}, \vec{a}_{j_2}, \ldots,\vec{a}_{j_s}$, 
and all $r$ vectors $\vec{a}_1,\vec{a}_2,\ldots, \vec{a}_r$ are orthogonal unit vectors, 
and stem from a complete basis of orthogonal unit vectors of ${\cal G}^1(I)$) 
precisely $k$ vectors of the two sets $\vec{a}_{j_1}, \vec{a}_{j_2}, \ldots,\vec{a}_{j_s}$, 
and $\vec{a}_1,\vec{a}_2,\ldots, \vec{a}_r$ must be in common. 
Otherwise the $m$-grade condition
for each specific $m$ is not fulfilled and the $\langle \rangle_m$-part is zero. For all non-zero $m$-grade
parts we can drop the $\langle \rangle_m$-bracket, use the reciprocal frame relationship of
P. \ref{P13}, and end up with $A_r=\vec{a}_1 \vec{a}_2 \ldots \vec{a}_r.$ 
How many such terms occur in the sum $\sum_{j_1<\ldots <j_s}\ldots$? This question can
be answered by considering that there are 
$
\left(
\begin{array}{c}
r \\
k
\end{array}
\right)
$
different choices of $k$ vectors from the set of $r$ vectors $\vec{a}_1,\vec{a}_2,\ldots, \vec{a}_r$. Then
there are $n-r$ vectors left from which we can still freely choose $s-k$ vectors. The freedom of
choice is reduced to $s-k$ vectors, because the set of $s$ vectors 
$\vec{a}_{j_1}, \vec{a}_{j_2}, \ldots,\vec{a}_{j_s}$ must have $k$ vectors in common with the $r$ vectors
 $\vec{a}_1,\vec{a}_2,\ldots, \vec{a}_r$, as just explained. There are 
$
\left(
\begin{array}{c}
n-r \\
s-k
\end{array}
\right)
$
different ways to choose $s-k$ vectors from $n-r$ vectors. 
\more
Hence 
\begin{align}
\sum_{j_1<\ldots <j_s} 
&\vec{a}^{j_s}\ldots\vec{a}^{j_2}\vec{a}^{j_1}
\langle 
   \vec{a}_{j_1}\vec{a}_{j_2}\ldots\vec{a}_{j_s}
   \vec{a}_1\vec{a}_2\ldots\vec{a}_r
\rangle_m
\nonumber \\
&= 
\left(
\begin{array}{c}
r \\
k
\end{array}
\right)
\left(
\begin{array}{c}
n-r \\
s-k
\end{array}
\right)
\delta^m_{r+s-2k}
A_{r}.
\label{eq:P49.D}
\end{align}
The result (\ref{eq:P49.D}) extends by linearity
 from the case of simple $A_r$ blades to general grade-homogeneous
multivectors $A_{\bar{r}}:$
\beqn
\partial_{\bar{s}}\langle X A_{\bar{r}}\rangle_m 
= \partial_{\bar{s}}\langle X_{\bar{s}} A_{\bar{r}}\rangle_m
=
\left(
\begin{array}{c}
r \\
k
\end{array}
\right)
\left(
\begin{array}{c}
n-r \\
s-k
\end{array}
\right)
\delta^m_{r+s-2k}
A_{\bar{r}}.
\label{eq:P49.3}
\end{equation}

The last part of our proof of P. \ref{P49} will be to show that the same result as in eq. (\ref{eq:P49.3}) 
holds true for $\langle A_{\bar{r}} \partial_{\bar{s}} \rangle_m X$. 
As before I again first assume $A_r$ to 
be a simple $r$-blade $A_r=\vec{a}_1 \vec{a}_2 \ldots \vec{a}_r$ of orthonormal vectors 
$\vec{a}_1,\vec{a}_2,\ldots, \vec{a}_r$. 
\newline
\more
We then have
\begin{align}
\langle A_{r} \partial_{\bar{s}} \rangle_m X
&\stackrel{\mbox{\scriptsize P. \ref{P14}}}{=}
\sum_J \left\langle \vec{a}_1 \vec{a}_2 \ldots \vec{a}_r \langle a_J \rangle_s \right\rangle_m
 \langle a^J \rangle_s \ast \partial_X \,X 
\nonumber \\
&\stackrel{\mbox{\scriptsize P. \ref{P45}}}{=}
\sum_J \left\langle \vec{a}_1 \vec{a}_2 \ldots \vec{a}_r \langle a_J \rangle_s \right\rangle_m
 \langle a^J \rangle_s
\nonumber \\
&= \sum_{j_1<\ldots <j_s} \langle \vec{a}_1 \vec{a}_2 \ldots \vec{a}_r
  \vec{a}_{j_1}\vec{a}_{j_2}\ldots\vec{a}_{j_s}  \rangle_m 
  \vec{a}^{j_s}\ldots\vec{a}^{j_2}\vec{a}^{j_1}
\nonumber \\
&=
\left(
\begin{array}{c}
r \\
k
\end{array}
\right)
\left(
\begin{array}{c}
n-r \\
s-k
\end{array}
\right)
\delta^m_{r+s-2k}
A_{r},
\label{eq:P49.4}
\end{align}
where I use in the last step the same argument as for the derivation of (\ref{eq:P49.D}), an important
part of which was the reciprocal frame relationship of P. \ref{P13}. Equation (\ref{eq:P49.4}) can again
be extended by linearity from simple $r$-blades $A_r$ to general grade-homogeneous
multivectors $A_{\bar{r}}$. This together with (\ref{eq:P49.3}) completes the full proof of P. {\ref{P49}}:
\beqn
\partial_{\bar{s}}\langle X A_{\bar{r}}\rangle_m
=
\langle A_{\bar{r}}\partial_{\bar{s}} \rangle_m X
=
\left(
\begin{array}{c}
r \\
k
\end{array}
\right)
\left(
\begin{array}{c}
n-r \\
s-k
\end{array}
\right)
\delta^m_{r+s-2k}
A_{\bar{r}}.
\end{equation}
\alles

\bpro
\normalfont
For $A=P(A)$, $K=\frac{1}{2}(r+s-|r-s|)$ and 
$\left(
\begin{array}{c}
i \\
j
\end{array}
\right) =0$ if $j>i$
\beqn
\partial_{\bar{s}} A_{\bar{r}} X_{\bar{s}}
= \sum_J \langle a^J \rangle_s A_{\bar{r}} \langle a_J \rangle_s
= \Gamma^r_s A_{\bar{r}},
\end{equation}
with 
\beqn
\Gamma^r_s = \sum_{k=0}^K (-1)^{rs-k} 
\left(
\begin{array}{c}
r \\
k
\end{array}
\right)
\left(
\begin{array}{c}
n-r \\
s-k
\end{array}
\right).
\end{equation}
\label{P50}
\epro
\prf{\ref{P50}}
\begin{align}
\langle A_{\bar{r}} X_{\bar{s}}\rangle_m
&\stackrel{\mbox{\scriptsize \cite{DH:CtoG}, p. 6, (1.20a)}}{=}
(-1)^{m(m-1)/2}\langle \tilde{A}_{\bar{r}} \tilde{X}_{\bar{s}}\rangle_m 
\nonumber \\
&=(-1)^{m(m-1)/2+r(r-1)/2+s(s-1)/2}\langle X_{\bar{s}} A_{\bar{r}} \rangle_m.
\label{eq:P49.1}
\end{align}
\newline
\more
Therefore
\begin{align}
A_{\bar{r}} X_{\bar{s}}
&\stackrel{\mbox{\scriptsize \cite{DH:CtoG}, p. 10, (1.36)}}{=}
\sum_{k=0}^K \langle A_{\bar{r}} X_{\bar{s}}\rangle_{|r-s|+2k}=
\sum_{k=0}^K \langle A_{\bar{r}} X_{\bar{s}}\rangle_{r+s-2k}
 \nonumber\\
&\stackrel{\mbox{\scriptsize (\ref{eq:P49.1})}}{=}
\sum_{k=0}^K \langle  X_{\bar{s}}A_{\bar{r}}\rangle_{r+s-2k}
\,\,(-1)^{
(r+s-2k)(r+s-2k-1)/2+r(r-1)/2+s(s-1)/2}.
\label{eq:P49.2}
\end{align}
The two different summations in line one of (\ref{eq:P49.2}) correspond to 
counting up as in \cite{DH:CtoG}, p. 10, (1.36):
\begin{align}
|r-s|, 
&|r-s|+2, |r-s|+4, \ldots 
\nonumber \\
&\ldots, |r-s|+2K=|r-s|+r+s-|r-s|=r+s,
\end{align}
and counting down as in \cite{DH:CtoG}, p. 58, after (2.38c):
\begin{align}
r+s, 
&r+s-2, r+s-4, \ldots 
\nonumber \\
&\ldots, r+s-2K=r+s-(r+s)+|r-s|= |r-s|.
\end{align}
\newline
\more
The exponent of $(-1)$ in (\ref{eq:P49.2}) can further be simplified to
\begin{align}
&(-1)^{(r+s-2k)(r+s-2k-1)/2+r(r-1)/2+s(s-1)/2} 
\nonumber \\
&=(-1)^{\frac{1}{2}(r^2+s^2+4k^2+2rs-4ks-4kr-r-s+2k+r^2-r+s^2-s)} 
\nonumber \\
&=(-1)^{\frac{1}{2}(2r^2+2s^2-2r-2s+2rs+2k)}
=(-1)^{(r^2+s^2-r-s+rs+k)}
\nonumber \\
&=(-1)^{r(r-1)+s(s-1)+rs-k+2k}
=(-1)^{rs-k}.
\end{align}
Equation (\ref{eq:P49.2}) simplifies therefore to
\begin{align}
A_{\bar{r}} X_{\bar{s}} = 
\sum_{k=0}^K \langle  X_{\bar{s}}A_{\bar{r}}\rangle_{r+s-2k}
\,\,(-1)^{rs-k}.
\label{eq:P50.1}
\end{align}
The $s$-grade derivative of (\ref{eq:P50.1}) gives
\begin{align}
\partial_{\bar{s}}A_{\bar{r}} X_{\bar{s}} 
&= 
\sum_{k=0}^K \partial_{\bar{s}}\langle  X_{\bar{s}}A_{\bar{r}}\rangle_{r+s-2k}
\,\,(-1)^{rs-k}
\nonumber \\
&\stackrel{\mbox{\scriptsize (\ref{eq:P49.3})}}{=}
\underbrace{
\sum_{k=0}^K
(-1)^{rs-k}
\left(
\begin{array}{c}
r \\
k
\end{array}
\right)
\left(
\begin{array}{c}
n-r \\
s-k
\end{array}
\right)}_{=\Gamma^r_s}
A_{\bar{r}} 
=\Gamma^r_s A_{\bar{r}}
\label{}
\end{align}
\newline
\more
To prove the first identity in P. \ref{P50} we write
\begin{align}
\partial_{\bar{s}}A_{\bar{r}} X_{\bar{s}} 
&\stackrel{\mbox{\scriptsize P. \ref{P14}}}{=}
\sum_J \langle a^J\rangle_s  
\underbrace{\langle a_J\rangle_s \ast \partial_{\bar{s}}}_{scalar}\,\,A_{\bar{r}} X_{\bar{s}} 
=
\sum_J \langle a^J\rangle_s A_{\bar{r}}  \left(\langle a_J\rangle_s \ast \partial_{\bar{s}} X_{\bar{s}}\right)
\nonumber \\
&\stackrel{\mbox{\scriptsize P. \ref{P36}}}{=}
\sum_J \langle a^J\rangle_s A_{\bar{r}} P_{
\stackrel{\mbox{\scriptsize $s$-dim.}}{\mbox{\scriptsize subspace}}}
(\langle a_J\rangle_s)
=\sum_J \langle a^J\rangle_s A_{\bar{r}}\langle a_J\rangle_s.
\end{align}
\alles

\bpro
\normalfont
\begin{equation}
\partial_{\bar{s}} A X = \partial A X_{\bar{s}}= \sum_{r=0}^n \Gamma^r_s A_{\bar{r}},
\end{equation}
with
$
\Gamma^r_s = \sum_{k=0}^K (-1)^{rs-k} 
\left(
\begin{array}{c}
r \\
k
\end{array}
\right)
\left(
\begin{array}{c}
n-r \\ 
s-k
\end{array}
\right)
$
and  $K=\frac{1}{2}(r+s-|r-s|).$
\label{P51}
\epro
\prf{\ref{P51}}
We will first proof that for an arbitrary but fixed grade $r$
\begin{equation}
\partial_{\bar{s}} A_{\bar{r}} X = \partial_{\bar{s}} A_{\bar{r}} X_{\bar{s}} = \partial A_{\bar{r}} X_{\bar{s}}.
\label{eq:P51.1}
\end{equation}
(For $r=0$ we end up with a scalar multiple of P. \ref{P47}.)
\begin{align}
\partial_{\bar{s}} A_{\bar{r}} X
&\stackrel{\mbox{\scriptsize P. \ref{P14}}}{=}
\sum_J \langle a^J\rangle_s 
\underbrace{\langle a_J\rangle_s \ast \partial_{\bar{s}}}_{scalar} A_{\bar{r}} X
=\sum_J \langle a^J\rangle_s A_{\bar{r}}
\underbrace{\langle a_J\rangle_s \ast \partial_{\bar{s}}X}  
\nonumber \\
&\stackrel{\mbox{\scriptsize P. \ref{P46}}}{=}
\sum_J \langle a^J\rangle_s A_{\bar{r}}\langle a_J\rangle_s
\stackrel{\mbox{\scriptsize P. \ref{P50}}}{=}
\partial_{\bar{s}} A_{\bar{r}} X_{\bar{s}}.
\end{align}
\newline
\more
The second identity in (\ref{eq:P51.1}) can be shown as follows:
\begin{align}
\partial A_{\bar{r}} X_{\bar{s}}
&\stackrel{\mbox{\scriptsize P. \ref{P14}}}{=}
\sum_{t=0}^n \sum_J \langle a^J\rangle_t 
\underbrace{\langle a_J\rangle_t \ast \partial}_{scalar} A_{\bar{r}} X_{\bar{s}}
\nonumber \\
&\stackrel{\mbox{\scriptsize \cite{DH:CtoG}, p. 13, (1.46)}}{=}
\sum_{t=0}^n \sum_J \langle a^J\rangle_t 
\underbrace{\langle a_J\rangle_t \ast \partial_{\bar{t}}}_{scalar} A_{\bar{r}} X_{\bar{s}}
\nonumber \\
&=\sum_{t=0}^n \sum_J \langle a^J\rangle_t A_{\bar{r}}
\underbrace{\langle a_J\rangle_t \ast \partial_{\bar{t}} X_{\bar{s}}}
\stackrel{\mbox{\scriptsize P. \ref{P3}}}{=}
\sum_t \sum_J \langle a^J \rangle_t A_{\bar{r}} \left\langle 
\langle a_J \rangle_t \ast \partial_{\bar{t}} X \right\rangle_s
\nonumber \\
&\stackrel{\mbox{\scriptsize P. \ref{P46}}}{=}
\sum_t \sum_J \langle a^J \rangle_t A_{\bar{r}} \left\langle P(\langle a_J \rangle_t) \right\rangle_s
=\sum_t \sum_J \langle a^J \rangle_t A_{\bar{r}} \underbrace{\left\langle \langle a_J \rangle_t \right\rangle_s}_{
=\delta_{ts}\langle a_J \rangle_t}
\nonumber \\
&=\sum_J \langle a^J \rangle_s A_{\bar{r}} \langle a_J \rangle_s
\stackrel{\mbox{\scriptsize P. \ref{P50}}}{=}
\partial_{\bar{s}} A_{\bar{r}} X_{\bar{s}}.
\end{align}
Performing the sum over all grades $r$ in (\ref{eq:P51.1}) we get
\begin{align}
\partial_{\bar{s}} A_{} X  = \partial A_{} X_{\bar{s}} = \sum_{r=0}^n \partial_{\bar{s}} A_{\bar{r}} X_{\bar{s}}
\stackrel{\mbox{\scriptsize P. \ref{P50}}}{=}
\sum_{r=0}^n \Gamma^r_s A_{\bar{r}},
\end{align}
\newline
\more
where I used the fact that according to Def. \ref{D15} and \cite{EH:VDcalc} (13):
\begin{equation}
A =  \sum_{r=0}^n A_{\bar{r}}.
\end{equation}
\alles

\bpro
\begin{align}
\partial_{\bar{s}}X \wedge A_{\bar{r}}=A_{\bar{r}}\wedge \partial_{\bar{s}}X = 
\left(
\begin{array}{c}
n-r \\ 
s
\end{array}
\right)
A_{\bar{r}}.
\end{align}
\label{P52}
\epro
\prf{\ref{P52}}
\begin{align}
\partial_{\bar{s}}X \wedge A_{\bar{r}}
&=
\partial_{\bar{s}} \left\langle X A_{\bar{r}} \right\rangle_{\mbox{\scriptsize max. grade}}
\stackrel{\mbox{\scriptsize (\ref{eq:P49.3})}}{=}
\partial_{\bar{s}} \left\langle X_{\bar{s}} A_{\bar{r}} \right\rangle_{\mbox{\scriptsize $r+s=$max. grade}}
\nonumber \\
&=\partial_{\bar{s}} \left\langle X A_{\bar{r}} \right\rangle_{\mbox{\scriptsize $r+s$}}
\stackrel{\mbox{\scriptsize P. \ref{P49}}}{=}
\left\langle A_{\bar{r}} \partial_{\bar{s}} \right\rangle_{\mbox{\scriptsize $r+s$}} X
= A_{\bar{r}} \wedge \partial_{\bar{s}} X,
\label{eq:P52.1}
\end{align}
\begin{align}
\partial_{\bar{s}}X \wedge A_{\bar{r}}
&\stackrel{\mbox{\scriptsize (\ref{eq:P52.1})}}{=}
\partial_{\bar{s}} \left\langle X A_{\bar{r}} \right\rangle_{\mbox{\scriptsize $r+s$}}
\stackrel{\mbox{\scriptsize P. \ref{P49}}}{=}
\left(
\begin{array}{c}
r \\
k
\end{array}
\right)
\left(
\begin{array}{c}
n-r \\
s-k
\end{array}
\right)
\delta^{m=r+s}_{r+s-2k}
A_{\bar{r}}
\nonumber \\
&\stackrel{\mbox{\scriptsize $k=0$}}{=}
\underbrace{\left(
\begin{array}{c}
r \\
0
\end{array}
\right)}_{=1}
\left(
\begin{array}{c}
n-r \\
s
\end{array}
\right)
A_{\bar{r}}
=
\left(
\begin{array}{c}
n-r \\
s
\end{array}
\right)
A_{\bar{r}}.
\end{align}
\alles

\bpro
\begin{align}
\partial_{\bar{s}} X \cdot A_{\bar{r}}
=A_{\bar{r}}\cdot \partial_{\bar{s}} X = 
\left\{
\begin{array}{c}
\left(
\begin{array}{c}
r \\
s
\end{array}
\right)A_{\bar{r}} \mbox{ if } 0<s\leq r
 \\
\left(
\begin{array}{c}
n-r \\
s-r
\end{array}
\right)A_{\bar{r}} \mbox{ if } 0<r\leq s
\end{array}
\right. .
\end{align}
\label{P53}
\epro
\more
\prf{\ref{P53}}
\begin{align}
\partial_{\bar{s}}X \cdot A_{\bar{r}}
&=
\partial_{\bar{s}} \left\langle X A_{\bar{r}} \right\rangle_{\mbox{\scriptsize min. grade}}
\stackrel{\mbox{\scriptsize (\ref{eq:P49.3})}}{=}
\partial_{\bar{s}} \left\langle X_{\bar{s}} A_{\bar{r}} \right\rangle_{\mbox{\scriptsize $|r-s|=$min. grade}}
\nonumber \\
&=\partial_{\bar{s}} \left\langle X A_{\bar{r}} \right\rangle_{\mbox{\scriptsize $|r-s|$}}
\stackrel{\mbox{\scriptsize P. \ref{P49}}}{=}
\left\langle A_{\bar{r}} \partial_{\bar{s}} \right\rangle_{\mbox{\scriptsize $|r-s|$}} X
= A_{\bar{r}} \cdot \partial_{\bar{s}}\,\, X
\label{eq:P53.1}
\end{align}
\begin{align}
\partial_{\bar{s}}X \cdot A_{\bar{r}}
&\stackrel{\mbox{\scriptsize (\ref{eq:P53.1})}}{=}
\partial_{\bar{s}} \left\langle X A_{\bar{r}} \right\rangle_{\mbox{\scriptsize $|r-s|$}}
\stackrel{\mbox{\scriptsize P. \ref{P49}}}{=}
\left(
\begin{array}{c}
r \\
k
\end{array}
\right)
\left(
\begin{array}{c}
n-r \\
s-k
\end{array}
\right)
\delta^{m=|r-s|}_{r+s-2k}
A_{\bar{r}}
\nonumber \\
&=\left\{
\begin{array}{c}
\left(
\begin{array}{c}
r \\
k
\end{array}
\right)
\left(
\begin{array}{c}
n-r \\
s-k
\end{array}
\right)A_{\bar{r}} 
\underbrace{\delta^{r-s}_{r+s-2k}}_{k=s}\mbox{ if } 0<s\leq r
 \\
\left(
\begin{array}{c}
r \\
k
\end{array}
\right)
\left(
\begin{array}{c}
n-r \\
s-k
\end{array}
\right)A_{\bar{r}} 
\underbrace{\delta^{s-r}_{r+s-2k}}_{k=r}\mbox{ if } 0<r\leq s
\end{array}
\right.
\nonumber \\
&=\left\{
\begin{array}{c}
\left(
\begin{array}{c}
r \\
s
\end{array}
\right)
\underbrace{\left(
\begin{array}{c}
n-r \\
s-s=0
\end{array}
\right)}_{=1}A_{\bar{r}} \mbox{ if } 0<s\leq r
 \\
\underbrace{\left(
\begin{array}{c}
r \\
r
\end{array}
\right)}_{=1}
\left(
\begin{array}{c}
n-r \\
s-k
\end{array}
\right)A_{\bar{r}} \mbox{ if } 0<r\leq s
\end{array}
\right. .
\end{align}
\alles

\bpro
\normalfont
For simple, $X$-independent $A_r$ the expansion
\begin{align}
\partial_{\bar{s}} X_{\bar{s}}
=&\partial_{\bar{s}} (X_{\bar{s}}A_r)A_r^{-1}
=\partial_{\bar{s}} X_{\bar{s}}\wedge A_r \,\,A_r^{-1}
+\partial_{\bar{s}} \langle X_{\bar{s}}\wedge A_r\rangle_{r+s-2} \,A_r^{-1}
\nonumber \\
&+\ldots+
\partial_{\bar{s}} \langle X_{\bar{s}}\wedge A_r\rangle_{|r-s|+2} \,A_r^{-1}
+\partial_{\bar{s}} X_{\bar{s}}\cdot A_r \,\,A_r^{-1}
\label{eq:P54.1}
\end{align}
is termwise equivalent to (assuming $0<s\leq r$)
\begin{align}
\left(
\begin{array}{c}
n \\
s
\end{array}
\right)
=
&\left(
\begin{array}{c}
r \\
0
\end{array}
\right)
\left(
\begin{array}{c}
n-r \\
s
\end{array}
\right)
+
\left(
\begin{array}{c}
r \\
1
\end{array}
\right)
\left(
\begin{array}{c}
n-r \\
s-1
\end{array}
\right)
\nonumber \\
&+\ldots+
\left(
\begin{array}{c}
r \\
s-1
\end{array}
\right)
\left(
\begin{array}{c}
n-r \\
1
\end{array}
\right)
+
\left(
\begin{array}{c}
r \\
s
\end{array}
\right)
\left(
\begin{array}{c}
n-r \\
0
\end{array}
\right).
\end{align}
\label{P54}
\epro
\prf{\ref{P54}}
 For simple $A_r$:
\begin{equation}
X_{\bar{s}}
\stackrel{\mbox{\scriptsize \cite{DH:CtoG}, p. 3, (1.3), \cite{EH:VDcalc}, (30)}}{=}
(X_{\bar{s}} A_r)A_r^{-1}.
\end{equation}
The expansion of $X_{\bar{s}} A_r$ is done according to \cite{DH:CtoG}, p. 10, (1.36).
\newline
\more
The correspondence of the first and last term are shown by P. \ref{P52} and P. \ref{P53}
respectively. In general each term in the right hand side expansion (\ref{eq:P54.1})
has the form
\begin{align}
\partial_{\bar{s}} \langle X_{\bar{s}} A_r \rangle_{r+s-2k} A_r^{-1}
&\stackrel{\mbox{\scriptsize P. \ref{P49}, (\ref{eq:P49.3})}}{=}
\left(
\begin{array}{c}
r \\
k
\end{array}
\right)
\left(
\begin{array}{c}
n-r \\
s-k
\end{array}
\right)
A_r A_r^{-1}
\nonumber \\
&=
\left(
\begin{array}{c}
r \\
k
\end{array}
\right)
\left(
\begin{array}{c}
n-r \\
s-k
\end{array}
\right),
\end{align}
with $k=0, \ldots, K$ and $K=\frac{1}{2}(r+s-|r-s|).$
The left hand side of (\ref{eq:P54.1}) is
\begin{equation}
\partial_{\bar{s}} X_{\bar{s}}
\stackrel{\mbox{\scriptsize P. \ref{P47}}}{=}
\left(
\begin{array}{c}
n \\
s
\end{array}
\right).
\end{equation}
The resulting binomial coefficient identity is known as theorem of addition: 
\cite{BS:TBdMa}, p. 105, (2.4).
\newline
\alles

\section{\cred Factorization}

Factorization relates functions of \textit{multivector} variables with corresponding functions of (several) 
\textit{lower grade}
content \textit{multivector} variables. In the simplest case the latter functions will
 just be functions of several \textit{vector} variables.

\bpro
\normalfont
For two multivector variables $A,B$:
\begin{align}
\partial_A G(A\wedge B)= \dot{\partial}_A(\dot{A}\wedge B)\ast \partial_U G_U(A \wedge B)
=\dot{\partial}_A\underline{G}(A \wedge B, \dot{A} \wedge B),
\end{align}
with 
\begin{equation}
\partial_U G_U (A \wedge B)
\stackrel{\mbox{\scriptsize P. \ref{P21}}}{=}
\partial_X G(X)|_{X=A \wedge B}
\end{equation}
\label{P55}
\epro
\prf{\ref{P55}}
For $f(A)=A\wedge B$ and $F=G(f(A))$:
\begin{equation}
\bar{f}(\partial_{X^{\prime}}) 
\stackrel{\mbox{\scriptsize P. \ref{P23}}}{=}
\dot{\partial}_A(\dot{A}\wedge B)\ast \partial_{X^{\prime}}.
\label{eq:P55.1}
\end{equation}
Therefore
\begin{align}
\partial_A F &= \partial_A G(f(A)) 
\stackrel{\mbox{\scriptsize P. \ref{P31}, P. \ref{P32}}}{=}
\bar{f}(\partial_{X^{\prime}}) 
G(X^{\prime})|_{X^{\prime}=f(A)=A\wedge B}
\nonumber \\
&\stackrel{\mbox{\scriptsize (\ref{eq:P55.1})}}{=}
\dot{\partial}_A(\dot{A}\wedge B)\ast \partial_{X^{\prime}}
G(X^{\prime})|_{X^{\prime}=A\wedge B}
\stackrel{\mbox{\scriptsize Def. \ref{D1}}}{=}
\dot{\partial}_A\underline{G}(X^{\prime},\dot{A}\wedge B )|_{X^{\prime}=A\wedge B}
\nonumber \\
&=\dot{\partial}_A\underline{G}(A\wedge B,\dot{A}\wedge B ).
\end{align}
\alles

\bpro
\normalfont
\begin{align}
\partial_B \partial_A G(A\wedge B)
=
&\,\,\stackrel{\circ}{\partial}_B \dot{\partial}_A (\dot{A}\wedge \stackrel{\circ}{B}) 
      \ast \partial_U \,\,G_U (A\wedge B)
\nonumber \\
&+ \stackrel{\circ}{\partial}_B ({A}\wedge \stackrel{\circ}{B})\ast \partial_V \,\, 
\dot{\partial}_A (\dot{A}\wedge{B})\ast \partial_U \,\,G_{UV} (A \wedge B),
\end{align}
with
\beqn
G_{UV}= V\ast \dot{\partial} \,\, U\ast \partial \,\,\dot{G}(A\wedge B).
\end{equation}
\label{P56}
\epro
\more
\prf{\ref{P56}}
\begin{align}
\partial_B \partial_A G(A\wedge B)
&\stackrel{\mbox{\scriptsize P. \ref{P55}}}{=}
\partial_B \left\{ \dot{\partial}_A (\dot{A} \wedge B)\ast \partial_U G_U (A\wedge B) \right\}
\nonumber \\
&\hspace{-2.3cm}\stackrel{\mbox{\scriptsize \cite{DH:CtoG}, p. 13 (1.44), P. \ref{P3}, Def. \ref{D12}(ii), P. \ref{P29}}}{=}
\,\,\stackrel{\circ}{\partial}_B \dot{\partial}_A (\dot{A}\wedge\stackrel{\circ}{B})
\ast \partial_U G_U(A\wedge B)
\nonumber \\
&\,\,\,\,\,\,\,\,\,\,+\stackrel{\circ}{\partial}_B \dot{\partial}_A (\dot{A}\wedge{B})
\ast  \stackrel{\circ}{(} \partial_U G_U(A\wedge B)). 
\label{eq:P56.1}
\end{align}

The second term on the right hand side of (\ref{eq:P56.1}) becomes
\begin{align}
\stackrel{\circ}{\partial}_B &\dot{\partial}_A (\dot{A}\wedge B)
\ast \stackrel{\circ}{(}\partial_U G_U(A\wedge B))
\stackrel{\mbox{\scriptsize P. \ref{P55}}}{=}
\nonumber \\
&{\color{black} \stackrel{\circ}{\partial}_B} 
{\color{black}\dot{\partial}_A (\dot{A} \wedge B) \ast }
\{
{\color{black}(A \wedge \stackrel{\circ}{B}) \ast }  \hspace{-0.9cm}
   \underbrace{{\color{black}\partial_{X^{\prime}}}
    {\color{black}\partial_U} G_U (X^{\prime})|_{X^{\prime}=A\wedge B}}_{
\stackrel{\mbox{\scriptsize Def. \ref{D12}(ii), P. \ref{P21}, Def. \ref{D10} }}{=}
\partial_V \partial_U G_{UV}(A\wedge B)
}
\hspace{-0.9cm}
\}
\nonumber \\
&=\stackrel{\circ}{\partial}_B\left\{ (A\wedge \stackrel{\circ}{B})\ast \partial_V \right\}
\,\,
\dot{\partial}_A \left\{ (\dot{A}\wedge B)\ast \partial_U  \right\} G_{UV} (A\wedge B).
\end{align}
\alles

\bpro
\normalfont
For a multivector function $G$ defined on ${\cal G}^{r}(I)$, i.e.
\begin{equation}
G(X) = G(\langle X \rangle_r),
\end{equation}
and for $A = \langle A\rangle_s$ in P. \ref{P54}, P. \ref{P55} we have $B=\langle B\rangle_{r-s}$ and
\begin{equation}
\partial_{A}G(A\wedge B) = B \cdot \partial_U G_U (A\wedge B).
\end{equation}
\label{P56a}
\epro
\more
\prf{\ref{P56a}}
\begin{align}
G(A\wedge B) = G(\langle A\wedge B\rangle_r)
\stackrel{\mbox{\scriptsize $G$ on ${\cal G}^{r}(I)$ }}{=}
 G(\langle \langle A\rangle_s \wedge B\rangle_r)
= G(\langle \langle A\rangle_s\langle B\rangle_{r-s}),
\label{eq:P56a.1}
\end{align}
and 
\begin{align}
\partial_A\underbrace{(A\wedge B)}_{=\langle A\wedge B\rangle_r}\ast 
\underbrace{\partial_X}_{=\langle \partial \rangle_r} 
&\stackrel{\mbox{\scriptsize \cite{DH:CtoG}, p. 13 (1.46)}}{=} 
\partial_A (\underbrace{A}_{=\langle A\rangle_s}\wedge \underbrace{B}_{=\langle B\rangle_{r-s}})\cdot 
\underbrace{\partial_X}_{=\langle \partial \rangle_r}
\nonumber \\
&\stackrel{\mbox{\scriptsize \cite{DH:CtoG}, p. 7 (1.25b)}}{=}
\partial_A (A\cdot (B\cdot \partial_X))
\nonumber \\
&\hspace{-0.7cm} \stackrel{\mbox{\scriptsize \cite{DH:CtoG}, p. 13 (1.46), $A = \langle A\rangle_s$}}{=} 
\partial_A (A\ast (B\cdot \partial_X))
\nonumber \\
&\hspace{0.55cm}\stackrel{\mbox{\scriptsize P. \ref{P46} }}{=}
B\cdot \partial_X.
\label{eq:P56a.2}
\end{align}
Using (\ref{eq:P56a.1}) and (\ref{eq:P56a.2}) we finally get
\begin{align}
\partial_A G(A\wedge B) 
&\stackrel{\mbox{\scriptsize P. \ref{P55} }}{=}
\dot{\partial}_A (\dot{A}\wedge B) \ast \partial_U G_U(A\wedge B)
\nonumber \\
&\hspace{-0.7cm}\stackrel{\mbox{\scriptsize similar to (\ref{eq:P56a.2}) }}{=}
B\cdot \partial_U G_U(A\wedge B).
\end{align}
\alles

\bpro
\normalfont
For a multivector function $G$ defined on ${\cal G}^{r}(I)$, i.e.
$$
G(X) = G(\langle X \rangle_r),
$$
and for $A = \langle A\rangle_s$ in P. \ref{P54}, P. \ref{P55} we have $B=\langle B\rangle_{r-s}$ and
\begin{align}
\partial_B \partial_A G(A\wedge B) = 
&\left(
\begin{array}{c}
r \\
r-s
\end{array}
\right)
\partial_U G_U (A\wedge B)
\nonumber \\
&+ (-1)^{s(r-s)}A\cdot \partial_V B\cdot \partial_U\,\, G_{UV}(A\wedge B).
\end{align}
\label{P57}
\epro
\prf{\ref{P57}}
The following three identities arise from proof \ref{P56a}:
\begin{align}
B=\langle B \rangle_{r-s},
\end{align}
\begin{align}
\dot{\partial}_A (\dot{A}\wedge B)\ast \partial_X = B\cdot \partial_X,
\label{eq:P57.2}
\end{align}
\begin{align}
\dot{\partial}_B(A\wedge \dot{B}) \ast \partial_X 
\stackrel{\mbox{\scriptsize \cite{EH:VDcalc} (40) }}{=}
\dot{\partial}_B(\dot{B}\wedge A) \ast \partial_X (-1)^{s(r-s)}
= A\cdot \partial_X (-1)^{s(r-s)}.
\label{eq:P57.3}
\end{align}
\more
We can therefore show that
\begin{align}
\stackrel{\circ}{\partial}_B \dot{\partial}_A (\dot{A}\wedge \stackrel{\circ}{B})\ast 
\partial_X 
\stackrel{\mbox{\scriptsize (\ref{eq:P57.2}) }}{=}
 \underbrace{\partial_B}_{\langle \partial_B\rangle_{r-s}} 
 \underbrace{B}_{\langle B \rangle_{r-s}}\cdot 
 \underbrace{\partial_X}_{\langle \partial_X \rangle_r} 
\stackrel{\mbox{\scriptsize P. \ref{P53} }}{=}
\left(
\begin{array}{c}
r \\
r-s
\end{array}
\right)
\partial_X.
\label{eq:P57.4}
\end{align}
We finally get for 
\begin{align}
\partial_B \partial_A G(A\wedge B)
\stackrel{\mbox{\scriptsize P. \ref{P56} }}{=}
&\stackrel{\circ}{\partial}_B \dot{\partial}_A (\dot{A}\wedge \stackrel{\circ}{B})\ast 
\partial_U G_U(A\wedge B)
\nonumber \\
&+\underbrace{\stackrel{\circ}{\partial}_B(A\wedge \stackrel{\circ}{B}) \ast {\partial}_V}_{
\mbox{\scriptsize similar to (\ref{eq:P57.3})}}
  \underbrace{\dot{\partial}_A (\dot{A}\wedge B) \ast \partial_U}_{
\mbox{\scriptsize similar to (\ref{eq:P57.2})}}
{G}_{UV}(A\wedge B)
\nonumber \\
\stackrel{\mbox{\scriptsize similar to (\ref{eq:P57.4}) }}{=}
&\left(
\begin{array}{c}
r \\
r-s
\end{array}
\right)
\partial_U G_U(A\wedge B)
\nonumber \\
&+ (-1)^{s(r-s)} A\cdot {\partial}_V   B\cdot \partial_U  \,\, {G}_{UV}(A\wedge B).
\end{align}
\alles

\bpro
\normalfont
For a linear multivector function $L=L(X)$ we have
\begin{align}
&\underline{L}
\stackrel{\mbox{\scriptsize Def. \ref{D1} }}{=}
L_U (X) = L(U),
 \\
&L_{UV}=0.
\end{align}
\label{P58}
\epro
\prf{\ref{P58}}
\begin{align}
\underline{L}
&\stackrel{\mbox{\scriptsize Def. \ref{D1} }}{=}
L_U(X) 
\stackrel{\mbox{\scriptsize Def. \ref{D1} }}{=}
\underbrace{U\ast \partial_X}_{scalar} L(X)
\nonumber \\
&\stackrel{\mbox{\scriptsize linearity }}{=}
L(U\ast \partial_X X )=L(P(U))
\stackrel{\mbox{\scriptsize for $P(U)=U$ }}{=}
L(U),
\label{eq:P58.3}
\end{align}
\begin{align}
L_{UV} 
&\stackrel{\mbox{\scriptsize Def. \ref{D10} }}{=}
V\ast \dot{\partial}_X U\ast \partial_X \dot{L}(X)
\stackrel{\mbox{\scriptsize (\ref{eq:P58.3}) }}{=}
V\ast \partial_XL(P(U))
\nonumber \\
&\stackrel{\mbox{\scriptsize linearity }}{=}
L(\underbrace{V\ast \partial_X P(U)}_{=0})
\stackrel{\mbox{\scriptsize linearity }}{=}0.
\end{align}
\alles

\bpro
\normalfont
A multivector function $L=L(U)$ is linear, if and only if
\beqn
\underline{L}=L(U).
\end{equation}
\label{P59}
\epro
\prf{\ref{P59}}
\\
$(\Rightarrow)$
\begin{align}
L=L(X) \mbox{ is linear } 
\stackrel{\mbox{\scriptsize P. \ref{P58} }}{\Rightarrow} 
\underline{L}=L(U),
\end{align}
$(\Leftarrow)$
\begin{align}
\underline{L}&=L(U) \Rightarrow L(\alpha U + \beta V)
= \underline{L}(X, \alpha U + \beta V) 
\nonumber \\
&\stackrel{\mbox{\scriptsize P. \ref{P4}, P. \ref{P5} }}{=}
\alpha \underline{L}(X,U) + \beta \underline{L}(X,V)
\nonumber \\
&\stackrel{\mbox{\scriptsize $\underline{L}=L(U)$ }}{=}
\alpha L(U) + \beta L(V)
\Rightarrow L \mbox{ is linear. }
\end{align}
\alles

\bpro[{\normalfont {\itshape factorization of derivative of lin. function}}]
\normalfont
If $G(A\wedge B)=L(A\wedge B)$ is a linear function on ${\cal G}(I),$ and $A=\langle A \rangle_s,$
$B=\langle B \rangle_{r-s},$
\beqn
\partial_B \partial_A L(A\wedge B)=
\left(
\begin{array}{c}
r \\
r-s
\end{array}
\right)
\partial_{\bar{r}}L.
\end{equation}
\label{P60}
\epro
\prf{\ref{P60}}
\begin{align}
\partial_B \partial_A L(A\wedge B)
&\stackrel{\mbox{\scriptsize P. \ref{P56}, P. \ref{P58} }}{=}
\stackrel{\circ}{\partial}_B \dot{\partial}_A (\dot{A}\wedge \stackrel{\circ}{B})
\ast \partial_U L_U(A\wedge B)
\nonumber \\
&\stackrel{\mbox{\scriptsize Def. \ref{D12}(ii), P. \ref{P21} }}{=}
\stackrel{\circ}{\partial}_B \dot{\partial}_A (\dot{A}\wedge \stackrel{\circ}{B})
\ast \partial_X L(X)|_{X=A\wedge B}
\nonumber \\
&\hspace{-2.3cm}\stackrel{\mbox{\scriptsize $A=\langle A \rangle_s,$ $B=\langle B \rangle_{r-s},$ 
   \cite{DH:CtoG}, p. 13 (1.45a) }}{=}
\stackrel{\circ}{\partial}_B \dot{\partial}_A \langle \dot{A}\wedge \stackrel{\circ}{B}\rangle_r
\ast \langle \partial_X\rangle_r L(X)|_{X=A\wedge B}
\nonumber \\
&\hspace{-1.1cm}\stackrel{\mbox{\scriptsize Proof \ref{P57}, P. \ref{P14}, Def. \ref{D15} }}{=}
\left(
\begin{array}{c}
r \\
r-s
\end{array}
\right)
\partial_{\bar{r}}L
= \left(
\begin{array}{c}
r \\
r-s
\end{array}
\right)
\partial L,
\end{align}
where the last equality holds if $L(X)= L(\langle X\rangle_r)$ on ${\cal G}^r(I).$
\newline
\alles

\bpro[{\normalfont {\itshape vector derivative factor}}]
\normalfont
For a linear function $L$:
\beqn
\partial_B \partial_{\vec{a}_1} L(\vec{a}_1\wedge B)= r \partial_{\bar{r}}L.
\end{equation}
\label{P61}
\epro
\prf{\ref{P61}}
\begin{align}
\partial_B \partial_{\vec{a}_1} L(\vec{a}_1\wedge B)
&\stackrel{\mbox{\scriptsize P. \ref{P60} }}{=}
\underbrace{\left(
\begin{array}{c}
r \\
r-(r-1)
\end{array}
\right)}_{=r}
\partial_{\bar{r}} L(\vec{a}_1\wedge B)
= r \partial_{\bar{r}} L.
\end{align}

\bpro
\normalfont
The function $\partial_{\vec{a}_1} L(\vec{a}_1\wedge B)$ is linear in $B$.
\label{P62}
\epro
\prf{\ref{P62}}
For $\alpha, \beta$ independent of $\vec{a}_1$:
\begin{align}
\partial_{\vec{a}_1} L(\vec{a}_1\wedge (\alpha B+ \beta C)) 
&\stackrel{\mbox{\scriptsize $L$ linear, \cite{EH:VDcalc} (23) }}{=}
\partial_{\vec{a}_1} \{\alpha L(\vec{a}_1\wedge B)+ \beta L(\vec{a}_1\wedge C) \}
\nonumber \\
&\stackrel{\mbox{\scriptsize \cite{EH:VDcalc} P. 51 }}{=}
\partial_{\vec{a}_1} \{\alpha L(\vec{a}_1\wedge B)\}
+\partial_{\vec{a}_1} \{\beta L(\vec{a}_1\wedge C)\}
\nonumber \\
&\hspace{-3cm}\stackrel{\mbox{\scriptsize \cite{EH:VDcalc} Def. 17(i) and eq. (21), (23), Def. \ref{D12}(i) }}{=}
\alpha \partial_{\vec{a}_1}  L(\vec{a}_1\wedge B)
+\beta \partial_{\vec{a}_1}  L(\vec{a}_1\wedge C).
\end{align}

\bpro
\normalfont
\begin{align}
\partial_r \ldots \partial_2 \partial_1 L(\vec{a}_1 \wedge \vec{a}_2 \ldots \wedge\vec{a}_r )
&= r! \,\,\partial_{\bar{r}} L 
\nonumber \\
&= \partial_r\wedge \ldots \wedge\partial_2 \wedge\partial_1 L(\vec{a}_1 \wedge \vec{a}_2 \ldots \wedge\vec{a}_r ),
\end{align}
with $\partial_k \equiv \partial_{\vec{a}_k}.$
\label{P63}
\epro
\prf{\ref{P63}}
\begin{align}
\partial_{\bar{r}}L
&\stackrel{\mbox{\scriptsize P. \ref{P61} }}{=}
\frac{1}{r}\,\, \partial_B \partial_1 L(\langle\vec{a}_1 \wedge 
  \underbrace{\vec{a}_2 \ldots \wedge\vec{a}_r}_{=B_{\overline{r-1}}} \rangle_r)
\nonumber \\
&= \frac{1}{r}\,\, \partial_B M(B) \,\,\,\mbox{ $M$ is linear according to P. \ref{P62} }
\nonumber \\
&=  \frac{1}{r}\,\, \partial_{\,\overline{r-1}}\, M
\stackrel{\mbox{\scriptsize P. \ref{P61} }}{=}
 \frac{1}{r}\,\,  \frac{1}{r-1}\,\, \partial_C \partial_2 M(\langle\vec{a}_2 \wedge 
  \underbrace{\vec{a}_3 \ldots \wedge\vec{a}_{r}}_{=C_{\overline{r-2}}} \rangle_{r-1})
\nonumber \\
&\stackrel{\mbox{\scriptsize P. \ref{P61} }}{=} 
\ldots
\stackrel{\mbox{\scriptsize P. \ref{P61} }}{=}
\frac{1}{r!}\partial_r Z (\vec{a}_r) \,\,\,\mbox{ $Z$ is linear according to P. \ref{P62}. }
\end{align}
\more
Hence
\begin{align}
r! \,\,\partial_{\bar{r}}\,\, L = \partial_r\,\, Z (\vec{a}_r)= \ldots = 
\partial_r \partial_{r-1} \ldots \partial_2 \partial_1 L(\vec{a}_1 \wedge \vec{a}_2 \ldots \wedge\vec{a}_r ).
\end{align}
We finally consider that
\begin{align}
\partial_r &\partial_{r-1} \ldots \partial_2 \partial_1 L(\vec{a}_1 \wedge \vec{a}_2 \ldots \wedge\vec{a}_r )
\nonumber \\
&=
\partial_r \partial_{r-1} \ldots (\partial_2\cdot  \partial_1+\partial_2 \wedge \partial_1) 
L(\vec{a}_1 \wedge \vec{a}_2 \ldots \wedge\vec{a}_r )
\nonumber \\
&\hspace{-0.65cm}\stackrel{\mbox{\scriptsize \cite{EH:VDcalc} Def. 17(i) }}{=}
\partial_r \partial_{r-1} \ldots \partial_2\cdot  \partial_1 
L(\vec{a}_1 \wedge \vec{a}_2 \ldots \wedge\vec{a}_r )
\nonumber \\
&\hspace{1.2cm}+
\partial_r \partial_{r-1} \ldots \partial_2 \wedge \partial_1 
L(\vec{a}_1 \wedge \vec{a}_2 \ldots \wedge\vec{a}_r )
\nonumber \\
&=
\partial_r \partial_{r-1} \ldots \partial_2 \wedge \partial_1 
L(\vec{a}_1 \wedge \vec{a}_2 \ldots \wedge\vec{a}_r ),
\label{eq:P63.3}
\end{align}
\more
because
\begin{align}
\partial_r &\partial_{r-1} \ldots \partial_2\cdot  \partial_1 
L(\vec{a}_1 \wedge \vec{a}_2 \ldots \wedge\vec{a}_r )
\nonumber \\
&=
\partial_r \partial_{r-1} \ldots 
\frac{1}{2}(\partial_2  \partial_1 +\partial_1  \partial_2)
L(\vec{a}_1 \wedge \vec{a}_2 \ldots \wedge\vec{a}_r )
\nonumber \\
&=
\partial_r \partial_{r-1} \ldots 
\frac{1}{2}(\partial_2  \partial_1 -
\underbrace{\partial_2  \partial_1}_{\stackrel{\mbox{\scriptsize relabeling}}{
   \mbox{\scriptsize $1 \leftrightarrow 2$}}}  
)
\hspace{0.5cm}L(\hspace{-0.4cm}\underbrace{\vec{a}_1 \wedge \vec{a}_2}_{
\stackrel{\mbox{\scriptsize interchanging}}{
   \mbox{\scriptsize $\vec{a}_1 \wedge \vec{a}_2=-\vec{a}_2 \wedge \vec{a}_1$}}    } 
\hspace{-0.4cm}\ldots \wedge\vec{a}_r )=0.
\label{eq:P63.4}
\end{align}
The factor $(-1)$ in the argument of $L$ caused by the interchange 
$\vec{a}_1 \wedge \vec{a}_2=-\vec{a}_2 \wedge \vec{a}_1$ and
the relabeling $(1 \leftrightarrow 2)$, can be factored out, because of the
linearity of $L.$
Applying the same consideration  to any pair of indizes $i,j \in \{1,2, \ldots,r \}$ of 
$$
\partial_r \partial_{r-1} \ldots \partial_2 \partial_1 L(\vec{a}_1 \wedge \vec{a}_2 \ldots \wedge\vec{a}_r )
$$
as to the pair $1,2$ in (\ref{eq:P63.3}) and (\ref{eq:P63.4}), 
we see that the vector derivatives on the left hand side are 
completely antisymmetric with respect to pairwise interchanges. Hence
\begin{align}
\partial_r \partial_{r-1} \ldots \partial_2 \partial_1 L(\vec{a}_1 \wedge \vec{a}_2 \ldots \wedge\vec{a}_r )
=\partial_r\wedge \ldots \wedge\partial_2 \wedge\partial_1 L(\vec{a}_1 \wedge \vec{a}_2 \ldots \wedge\vec{a}_r ).
\end{align}
\alles

\section{\cred Simplicial Derivatives and Variables}

\bdf[{\normalfont {\itshape simplicial variable, simplicial derivative}}]
\normalfont
A \textit{simplicial variable} $a_{(r)}$ is a simple blade of the form
\beqn
a_{(r)}\equiv \vec{a}_1  \wedge\ldots \wedge \vec{a}_r.
\end{equation}
The \textit{simplicial derivative} $\partial_{(r)}$ is defined as
\beqn
\partial_{(r)}\equiv \frac{1}{r!} \partial_{\vec{a}_1}  \wedge\ldots \wedge \partial_{\vec{a}_r}.
\end{equation}
\label{D64}
\edf

\bpro[{\normalfont {\itshape equivalence}}]
\normalfont
For linear functions $L$ the $r$-vector derivative (right hand side) is equivalent to the simplicial
derivative (left hand side):
\begin{equation}
\partial_{(r)}L(a_{(r)})= \partial_{\bar{r}} L.
\end{equation}
\label{P65}
\epro
\prf{\ref{P65}}
P. \ref{P63} and Def. \ref{D64}.

\bpro
\beqn
\partial_{(r)}a_{(r)}= \partial_{\bar{r}} X_{\bar{r}}=
\left(
\begin{array}{c}
n \\
r
\end{array}
\right).
\end{equation}
\label{P66}
\epro
\prf{\ref{P66}}
\begin{equation}
\partial_{(r)}a_{(r)}
\stackrel{\mbox{\scriptsize P. \ref{P65} }}{=}
\partial_{\bar{r}} X
\stackrel{\mbox{\scriptsize P. \ref{P47} }}{=} 
\partial_{\bar{r}} X_{\bar{r}}
\stackrel{\mbox{\scriptsize P. \ref{P47} }}{=}  
\left(
\begin{array}{c}
n \\
r
\end{array}
\right).
\end{equation}

\bpro
\beqn
\partial_{(r)}\left(a_{(r)}\right)^2= (r+1)a_{(r)}.
\end{equation}
\label{P67}
\epro
\brem
\normalfont
Though I lack the general proof so far, the expamles of one, two and three dimensions
are quite instructive.
\\
1-D
\\
\begin{equation}
\partial_{(1)}\left(a_{(1)}\right)^2=\partial_{\vec{a}_1}(\vec{a}_1)^2
\stackrel{\mbox{\scriptsize P. \ref{P47}, \cite{EH:VDcalc} P. 69 }}{=}
2 \vec{a}_1.
\end{equation}
\\
2-D
\\
\begin{align}
\partial_{(2)}\left(a_{(2)}\right)^2
&=\frac{1}{2} \partial_2 \wedge \partial_1 (\vec{a}_1\wedge \vec{a}_2)^2
=\frac{1}{2} \partial_2 \wedge \partial_1 [(\vec{a}_1\cdot \vec{a}_2)^2-\vec{a}_1^2\vec{a}_2^2]
\nonumber \\
&=\frac{1}{2} \partial_2 \wedge \partial_1 (\vec{a}_1\cdot \vec{a}_2)^2
  -\frac{1}{2}(\partial_2\vec{a}_2^2)\wedge(\partial_1\vec{a}_1^2)
\nonumber \\
&=\frac{1}{2}\partial_2\wedge \{2(\vec{a}_1\cdot \vec{a}_2)\partial_1(\vec{a}_1\cdot \vec{a}_2)  \}
  -\frac{1}{2}2\vec{a}_2\wedge(2\vec{a}_1)
\nonumber \\
&=\partial_2(\vec{a}_1\cdot \vec{a}_2)\wedge \vec{a}_2 + 2 \vec{a}_1 \wedge \vec{a}_2
\nonumber \\
&=\vec{a}_1 \wedge \vec{a}_2 +\vec{a}_1\cdot \vec{a}_2 
  \underbrace{\partial_2\wedge\vec{a}_2}_{
\stackrel{\mbox{\scriptsize $=0$}}{\mbox{\scriptsize \cite{EH:VDcalc}, P. 67}} }
  +2\vec{a}_1 \wedge \vec{a}_2
=3\vec{a}_1 \wedge \vec{a}_2.
\end{align}
\more
3-D
\\
\begin{align}
&\partial_{(3)}\left(a_{(3)}\right)^2
=\frac{1}{3!} \partial_3 \wedge \partial_2 \wedge \partial_1 (\vec{a}_1\wedge \vec{a}_2\wedge \vec{a}_3)^2
\nonumber \\
&=\frac{1}{3!} \partial_3 \wedge \partial_2 \wedge \partial_1 
(-\vec{a}_1^2\vec{a}_2^2\vec{a}_3^2
  +\vec{a}_1^2(\vec{a}_2\cdot \vec{a}_3)^2
  +\vec{a}_2^2(\vec{a}_3\cdot \vec{a}_1)^2
  +\vec{a}_3^2(\vec{a}_1\cdot \vec{a}_2)^2
  \nonumber \\
  &\hspace{0.5cm}-2\vec{a}_1\cdot \vec{a}_2 \vec{a}_2\cdot \vec{a}_3 \vec{a}_3\cdot \vec{a}_1 )
\nonumber \\
&=\frac{1}{6} (-(\partial_3\vec{a}_3^2) \wedge (\partial_2\vec{a}_2^2) \wedge (\partial_1\vec{a}_1^2)
  +(\partial_1\vec{a}_1^2)\wedge (\partial_3 \wedge \partial_2)(\vec{a}_2\cdot \vec{a}_3)^2
  \nonumber \\
  &\hspace{0.5cm}-(\partial_2\vec{a}_2^2)\wedge (\partial_3 \wedge \partial_1)(\vec{a}_1\cdot \vec{a}_3)^2
  +(\partial_3\vec{a}_3^2)\wedge (\partial_2 \wedge \partial_1)(\vec{a}_1\cdot \vec{a}_2)^2
  \nonumber \\
  &\hspace{0.5cm}-2(\partial_3 \wedge \partial_2 \wedge \partial_1) 
     [\vec{a}_1\cdot \vec{a}_2 \vec{a}_2\cdot \vec{a}_3 \vec{a}_3\cdot \vec{a}_1]
              )
\nonumber \\
&\stackrel{\mbox{\scriptsize 2-D}}{=}
\frac{1}{6}\{
   -8 \vec{a}_3\wedge \vec{a}_2\wedge \vec{a}_1
   +2\vec{a}_1\wedge[2(\vec{a}_2\wedge \vec{a}_3)]
   -2\vec{a}_2\wedge[2(\vec{a}_1\wedge \vec{a}_3)]
   \nonumber \\
   &\hspace{0.5cm}+2\vec{a}_3\wedge[2(\vec{a}_1\wedge \vec{a}_2)]
   -2\partial_3 \wedge \partial_2 \wedge(\vec{a}_2 \vec{a}_2\cdot \vec{a}_3 \vec{a}_3\cdot \vec{a}_1
                                         + \vec{a}_3 \vec{a}_1\cdot \vec{a}_2 \vec{a}_2\cdot \vec{a}_3)
           \}
\nonumber \\
&\stackrel{\mbox{\scriptsize \cite{EH:VDcalc}, P. 67}}{=}
\frac{1}{6}\{
   8 \vec{a}_1\wedge \vec{a}_2\wedge \vec{a}_3
   +(4+4+4)\vec{a}_1\wedge \vec{a}_2\wedge \vec{a}_3
   \nonumber \\
   &\hspace{1cm}-2\partial_3 \wedge (\vec{a}_3\wedge \vec{a}_2(\vec{a}_3\cdot \vec{a}_1)
                        +\vec{a}_1\wedge \vec{a}_3(\vec{a}_2\cdot \vec{a}_3)
                        +\underbrace{\vec{a}_3\wedge \vec{a}_3}_{=0}(\vec{a}_1\cdot \vec{a}_2) )
           \}
\nonumber \\
&\stackrel{\mbox{\scriptsize \cite{EH:VDcalc}, P. 67}}{=}
\frac{1}{6}\{
   20 \vec{a}_1\wedge \vec{a}_2\wedge \vec{a}_3  
  -2(\vec{a}_1\wedge \vec{a}_3\wedge \vec{a}_2 + \vec{a}_2\wedge \vec{a}_1\wedge \vec{a}_3)
\nonumber \\
&=\frac{24}{6} \vec{a}_1\wedge \vec{a}_2\wedge \vec{a}_3 
= 4 \vec{a}_1\wedge \vec{a}_2\wedge \vec{a}_3.\hspace{1cm} \mbox{\alles}
\end{align}
\label{R67}
\erem

\bpro
\normalfont
\begin{align}
\left(
\begin{array}{c}
r+s \\
s
\end{array}
\right)
\partial_{(r+s)L(a_{(r+s)})}
= \partial_{(r)}\wedge \partial_{(s)} L(a_{(r)}\wedge a_{(s)})
= \partial_{(r)} \partial_{(s)} L(a_{(r)}\wedge a_{(s)})
\end{align}
for $\partial_{(r)}$ and $a_{(r)}$ as in Def. \ref{D64} and $L$ a linear function.
\label{P68}
\epro
\prf{\ref{P68}}
\begin{align}
a_{(r+s)}
\stackrel{\mbox{\scriptsize Def. \ref{D64}}}{=} 
a_{(r)}\wedge a_{(s)},
\end{align}
\begin{align}
\partial_{(r+s)}
&\stackrel{\mbox{\scriptsize Def. \ref{D64}}}{=} 
\frac{1}{(r+s)!} \partial_{r+s}\wedge \ldots \wedge \partial_{s+1}\wedge \partial_{s}\wedge \ldots \wedge \partial_1
\nonumber \\
&=\frac{r!s!}{(r+s)!} (\frac{1}{r!}\partial_{r+s}\wedge \ldots \wedge \partial_{s+1})\wedge
                     (\frac{1}{s!}\partial_{s}\wedge \ldots \wedge \partial_1)
\nonumber \\
&=\left( \frac{(r+s)!}{(r+s-s)!s!}\right)^{-1}
   \partial_{(r)}\wedge \partial_{(s)} 
=
\left(
\begin{array}{c}
r+s \\
s
\end{array}
\right)^{-1}
\partial_{(r)}\wedge \partial_{(s)}.
\end{align}
\more
Hence
\begin{align}
\left(
\begin{array}{c}
r+s \\
s
\end{array}
\right)
\partial_{(r+s)}=\partial_{(r)}\wedge \partial_{(s)},
\end{align}
and 
\begin{align}
\left(
\begin{array}{c}
r+s \\
s
\end{array}
\right)
\partial_{(r+s)}L(a_{(r+s)}) &=\partial_{(r)}\wedge \partial_{(s)}L(a_{(r)}\wedge a_{(s)})
\nonumber \\
&=\partial_{(r)}\partial_{(s)}L(a_{(r)}\wedge a_{(s)}).
\end{align}
The last equality holds because of the symmetry of the argument $a_{(r)}\wedge a_{(s)}$.
\newline
\alles

\bpro
\normalfont
For linear functions $L$:
\begin{align}
\partial_B \partial_A L(A\wedge B)= \partial_B \wedge \partial_A L(A\wedge B).
\end{align}
\label{P69}
\epro
\prf{\ref{P69}}
\begin{align}
\partial_B \partial_A L(A\wedge B)
&\stackrel{\mbox{\scriptsize P. \ref{P60}}}{=}
\left(
   \begin{array}{c}
      r \\
      r-s
   \end{array}
\right)
\partial_{\bar{r}}L
\stackrel{\mbox{\scriptsize P. \ref{P65}}}{=}
\left(
   \begin{array}{c}
      r \\
      r-s
   \end{array}
\right)
\partial_{(r)}L(a_{(r)})
\nonumber \\
&\stackrel{\mbox{\scriptsize $
\left(
   \begin{array}{c}
      r \\
      r-s
   \end{array}
\right)
=
\left(
   \begin{array}{c}
      r \\
      s
   \end{array}
\right)
$}}{=}
\left(
   \begin{array}{c}
      r-s+s \\
      s
   \end{array}
\right)
\partial_{(r-s+s)}L(a_{(r-s+s)})
\nonumber \\
&\stackrel{\mbox{\scriptsize P. \ref{P68}}}{=}
\partial_{(r-s)}\wedge \partial_{(s)} L(a_{(s)}\wedge a_{(r-s)})
 \\
&\stackrel{\mbox{\scriptsize P. \ref{P65}}}{=}
\partial_{\overline{r-s}}\wedge \partial_{(s)} L(a_{(s)}\wedge B_{\overline{r-s}})
\nonumber \\
&\stackrel{\mbox{\scriptsize P. \ref{P65}}}{=}
\partial_{\overline{r-s}}\wedge \partial_{\overline{s}} L(A_{\overline{s}}\wedge B_{\overline{r-s}})
=
\partial_B\wedge \partial_A L(A_{\overline{s}}\wedge B_{\overline{r-s}}).
\nonumber
\end{align}

\bpro
\normalfont
For a linear function $F_r=F_r(\vec{a}_1,\vec{a}_2,\ldots,\vec{a}_r )$ of $r$ vector variables
\begin{align}
\vec{a}_k\cdot \partial_{\vec{a}_k}F_r=F_r, \hspace{1cm}  (k=1,\ldots, r).
\end{align}
\label{P70}
\epro
\prf{\ref{P70}}
\begin{align}
\vec{a}_k\cdot \partial_{\vec{a}_k}F_r
&=
\vec{a}_k\cdot \partial_{\vec{a}_k}F_r(\vec{a}_1,\vec{a}_2,\ldots,\vec{a}_k, \ldots,\vec{a}_r)
\nonumber \\
&\hspace{-1.4cm}\stackrel{\mbox{\scriptsize P. \ref{P58}, P. \ref{P59}, $X\equiv\vec{a}_k=U $}}{=}
F_r(\vec{a}_1,\vec{a}_2,\ldots,\vec{a}_k, \ldots,\vec{a}_r),
\end{align}
for all $k=1,\ldots, r.$

\bdf[{\normalfont {\itshape skew-symmetrizer}}]
\normalfont
\begin{align}
G_r\equiv\frac{1}{r!}(\vec{a}_1\wedge \ldots \wedge\vec{a}_r)\cdot (\partial_r\wedge \ldots \wedge \partial_1)F_r
= a_{(r)} \cdot \partial_{(r)} F_r,
\end{align}
for  a linear function $F_r$ of $r$ vector variables.
\label{D70}
\edf

\bpro
\normalfont
\begin{align}
G_r=G_r(\vec{a}_1,\ldots,\vec{a}_r)
\end{align}
is according to Def. \ref{D70} skew-symmetric in its arguments. If $F_r$ is skew-symmetric, then
\begin{align}
G_r=F_r.
\end{align}
\label{P71}
\epro
\more
\prf{\ref{P71}}
\begin{align}
&(\vec{a}_1\wedge \ldots \wedge\vec{a}_r)\cdot (\partial_r\wedge \ldots \wedge \partial_1)F_r
\nonumber \\
&\hspace{-0cm}\stackrel{\mbox{\scriptsize \cite{DH:CtoG} p. 7 (1.25b), p. 11 (1.38)} }{=}
(\vec{a}_1\wedge \ldots \wedge\vec{a}_{r-1})\cdot [\vec{a}_r\cdot(\partial_r\wedge \ldots \wedge \partial_1)]F_r
\nonumber \\
&\hspace{-0cm}\stackrel{\mbox{\scriptsize \cite{DH:CtoG} p. 11 (1.38)} }{=}
\sum_{k=1}^{r}(-1)^{r-k} (\vec{a}_1\wedge \ldots \wedge\vec{a}_{r-1})\cdot 
   (\partial_r\wedge \ldots\wedge\check{\partial}_k\wedge\ldots \wedge \partial_1)
    \nonumber \\
    &\hspace{40mm}\vec{a}_r\cdot\partial_k F_r(\vec{x}_1,\ldots,\vec{x}_k,\ldots,\vec{x}_r)
\nonumber \\
&\stackrel{\mbox{\scriptsize $F$ linear, P. \ref{P70}} }{=}      
\sum_{k=1}^{r}(-1)^{r-k} (\vec{a}_1\wedge \ldots \wedge\vec{a}_{r-1})\cdot 
   (\partial_r\wedge \ldots\wedge\check{\partial}_k\wedge\ldots \wedge \partial_1)
    \nonumber \\
    &\hspace{40mm}
      F_r(\vec{x}_1,\ldots,\hspace{-0.2cm}\stackrel{\mbox{\scriptsize position $k$}}{\vec{a}_r}\hspace{-0.2cm}, 
           \ldots,\vec{x}_r)
\nonumber \\
&\stackrel{\mbox{\scriptsize $F$ skew-symmetric} }{=}
\sum_{k=1}^{r}
   (\vec{a}_1\wedge \ldots \wedge\vec{a}_{r-1})\cdot 
   (\partial_r\wedge \ldots\wedge\check{\partial}_k\wedge\ldots \wedge \partial_1)
    \nonumber \\
    &\hspace{40mm}
   F_r(\vec{x}_1,\ldots,\check{\vec{x}}_k,\ldots,\vec{x}_r,\vec{a}_r)
\nonumber \\
&\stackrel{\stackrel{\mbox{\scriptsize relabeling of}}{\mbox{\scriptsize $\vec{x}_1,\ldots,\vec{x}_r$} }}{=}
r (\vec{a}_1\wedge \ldots \wedge\vec{a}_{r-1})\cdot 
  (\partial_{r-1}\wedge \ldots \wedge \partial_1)
   F_r(\vec{x}_1,\ldots,\vec{x}_{r-1},\vec{a}_r)
\nonumber \\
&\stackrel{\mbox{\scriptsize P. \ref{P70}} }{=}
   r(\vec{a}_1\wedge \ldots \wedge\vec{a}_{r-1})\cdot (\partial_{r-1}\wedge \ldots \wedge \partial_1)
         \vec{a}_r \cdot \partial_r F_r(\vec{x}_1,\ldots,\vec{x}_{r-1},\vec{x}_r)
\nonumber \\
&=\ldots =
r! \vec{a}_1 \cdot \partial_1 \ldots \vec{a}_r \cdot \partial_r F_r
\stackrel{\mbox{\scriptsize P. \ref{P70}} }{=}
r! F_r,
\end{align}
\more
 where the symbols
$\check{\vec{x}}_k,$ $\check{\partial}_k$ in the above mean that ${\vec{x}}_k$ and
$\partial_k=\partial_{\vec{x}_k}$ are to be left out.
\alles

\bpro[{\normalfont {\itshape canonical form of alternating linear form}}]
\normalfont
The simplicial derivative of an alternating linear $r$-form $\alpha_r=\alpha_r(\vec{a}_1,\ldots,\vec{a}_{r})$ 
(\cite{DH:CtoG}, pp. 33 ff.)
yields a $r$-vector
\begin{equation}
A^+_r\equiv  \partial_{(r)}\alpha_r,
\end{equation}
with the property
\begin{equation}
\alpha_r = (\vec{a}_1\wedge\ldots\wedge\vec{a}_{r})\cdot A^+_r 
= a_{(r)}\cdot A^+_r,
\label{eq:P72.2}
\end{equation}
 also called the \textit{canonical form} of an alternating linear $r$-form $\alpha_r.$
\label{P72}
\epro
\prf{\ref{P72}}
Def. \ref{D64} and Def. \ref{D12}(i) (for vector derivatives) show that the above defined $A_r^+$ must be 
an $r$-vector. 
\begin{align}
a_{(r)}\cdot A_r^+ = a_{(r)}\cdot \partial_{(r)}\alpha_r
= \alpha_r.
\end{align}
The last equality is due to P. \ref{P70}, and the linearity and skew-symmetry of $\alpha_r.$
\newline
\alles

\section{\cred Conclusion}

In order to make the treatment in the future more self-contained, it may be useful to 
compile a set of important geometric algebra 
relationships, which are necessary for multivector differential calculus, and often referred to in
this paper.  

It may be of interest to notice that there is a MAPLE package software implementation of the multivector derivative
and multivector differential~\cite{Cam:MAPLE}.

After discussing vector differential calculus \cite{EH:VDcalc}  and multivector differential calculus
in some detail, I would like to proceed and produce a similar discussion on directed integration
of multivector functions in the future.

\subsection*{Acknowledgements}
I first of all thank God for the joy of studying his creation: 
\begin{quotation}
{\color{black}"\ldots since the creation of the world God's invisible qualities - his eternal power and divine nature - 
have been clearly seen, being understood from what has been made \ldots"}\cite{P:Rom}.
\end{quotation}

{\cred I thank my wife for 
encouragement,}  H. Ishi for corrections and improvements, G. Sobczyk for some hints,
 and O. Giering and J.S.R. Chisholm for attracting me 
to geometry. Fukui University provided a good research environment.
I thank K. Shinoda (Kyoto) for his prayerful support.


\end{document}

%% file: Makro98.tex
\newcommand{\cl}{C \kern -0.1em \ell}     

\newcommand{\ut}[1]{{\setbox0=\hbox{$#1$}\mathsurround=0pt
       \rlap{\raisebox{-0.8\dp0}{\raisebox{-0.8ex}
       {\kern -0.15ex\hbox{$\tiny\sim$}\kern 0.15ex}}}#1}}
\newcommand{\uti}[1]{{\setbox0=\hbox{$#1$}\mathsurround=0pt
       \rlap{\raisebox{-0.8\dp0}{\raisebox{-0.8ex}
       {\kern -0.3ex\hbox{$\tiny\sim$}\kern 0.3ex}}}#1}}

\newdimen\arrayruleHwidth                 
     \makeatletter
     \def\Hline{\noalign{\ifnum0=`}\fi\hrule \@height \arrayruleHwidth
         \futurelet \@tempa\@xhline}
     \makeatother